\documentclass[a4paper,10pt]{amsart}

\usepackage{amssymb}

\newcommand{\B}{\mathcal{B}}

\newcommand{\ram}{r_{T}'}
\newcommand{\vsp}{\vspace*{.5cm}}

\newcommand{\Bm}{\mathcal{B}_{m}}

\newcommand{\Top}{\mathcal{T}}
\begin{document}

\title{SOME REMARKS ON PREDICTION MODELS}
\author{T. Dahn}
\footnote{Lund University}
\maketitle

\section{Introduction}
This article will deal with global models, particularly for the purpose of prediction of a geometric event. It contains a comparison of some different global models and problems known to us. We will localize the symbol, using (analytic) continuation, to a given geometric event and we discuss a transmission property, that will transform the localization to a normal model. The localization is an inverse lifting principle,
$F_{T}(\gamma)(\zeta) \rightarrow \gamma_{T}(\zeta) \rightarrow \zeta_{T}$, where given a geometric ideal, we look for a set of systems of trajectories $\gamma_{T}(\zeta)$, for which the domain include the given event. Further, we discuss the events that can be reached using a set of fixed systems. Particularly, we consider the situation when the systems are considered at different points in time. 
The characteristic properties in the representations we are discussing are the approximation property, transmission property and the interpolation property. We assume always in this paper the original symbol has summable phase, $\log \mid f \mid \in L^{1}$.

\section{Topological considerations}

For reference we consider the spaces of monotropic functionals. We will use the notation $f \sim_{m}0$  explained as follows.
Between the spaces $\dot{\B}(\mathbb{R}^{n})$ and
$\B(\mathbb{R}^{n})$, we consider over an $\epsilon-$ neighborhood
of the real space, the space $\Bm$ of $C^{\infty}-$ functions
bounded in the real infinity by a small constant with all
derivatives. Thus, consider $D^{\alpha}\phi - \mu_{\alpha}
\rightarrow 0$ in the real infinity, for all $\alpha$ and
$\mu_{\alpha}$ constants. Obviously, the space of monotropic functionals $\Bm' \subset
\mathcal{D}_{L^{1}}'$, why $T \in \Bm'$ has representation
$\sum_{\mid \alpha \mid \geq k} D^{\alpha}f_{\alpha}$ with
$f_{\alpha} \in L^{1}$. If $T \in \mathcal{D}_{L^{1}}'$ and $\phi
\in \Bm$, there is a $S \in \Bm'$ such that $S=T$ over $\Bm$. We
have that $\mathbb{R}^{n}=\cup^{\infty}_{j=0} K_{j}$, for compact
sets $K_{j}$. Let $\Phi_{j,1}=(S-T) |_{K_{j}} \in \mathcal{E}'
\subset \Bm'$ and $\Phi_{j,2}=\Phi-\Phi_{j,1}$. We chose $S$ such
that $\Phi_{j,1}=0$ for all $j$ and $\Phi_{j,2} \in \Bm'$. This
gives existence of a functional $S$ such that $S(\phi)=\sum_{\alpha}
f_{\alpha}(x)dx=\lim_{j \rightarrow \infty} T_{j}(\phi)$, where the
limit is taken in $\mathcal{D}_{L^{1}}'$ (\cite{Schwartz66})

\vsp

Consider in particular the example with an algebraic homomorphy $h$, such that $h^{*}=h$
and $h^{2}=1$, $f=e^{\phi}$.
More precisely, if $h(\frac{f}{g})=0$, we have $\mid f \mid=\mid g \mid$, that is 
$h(\frac{f}{g}) \rightarrow 0$, then $\mid \mid f \mid - \mid g \mid \mid \rightarrow 0$ and 
$d h(\frac{f}{g}) \rightarrow 0$, then $\mid f - g \mid \rightarrow 0$, this close to the infinity. 
Further, if we assume $\widehat{h}=H$, then $H(\frac{\phi}{\psi}) \rightarrow 0$, then 
$\mid \phi \mid = \mid \psi \mid$ and $d H(\frac{\phi}{\psi})=0$ implies $\mid \phi - \psi \mid \rightarrow 0$. 
Further, $\widehat{d h}(\frac{\phi}{\psi}) \rightarrow 0$, means $(d h)(e^{\frac{\phi}{\psi}}) \rightarrow 0$, 
that is $\frac{\phi}{\psi} \rightarrow 0$, in the infinity.

\subsubsection*{The summable distributions}
Let $\mathcal{B}_{\alpha}=\{ f \quad \mid \xi \mid^{\alpha} \mid \widehat{f} \mid \in L^{1} \}$ (\cite{Nilsson72}). Assume $(L_{f})=\{ N \quad Nf \in \mathcal{B}_{0} \quad \{ \widehat{Nf} =0 \} \subset \{ \widehat{f}=0 \} \}$. We assume that if $N_{1},N_{2} \in (L_{f})$, then $N_{1}N_{2}=N_{2}N_{1} \in (L_{f})$. Further, if $\tilde{f}_{0}=N_{0}f$, then $\tilde{f}_{1}=N_{1} f$, then we have existence of $N_{2} \in (L_{f})$, such that $N_{2}N_{0}f=N_{1}f$. For instance if $\parallel \eta \parallel^{*}$ is a norm (locally 1-1), then $\parallel \eta_{0}+\eta_{2} \parallel^{*}=\parallel \eta_{1} \parallel^{*}$ and $\parallel \eta_{0}+\eta_{2}-\eta_{1} \parallel^{*}=0$, that is the interpolation property is interpreted as a unique $\eta$ such that $\parallel \eta \parallel^{*}=0$. However if we only have that the dual norm is $0$, we can only conclude that the type is $0$, that is the functional is of real type, the continuation is through translation. We note that $\parallel \eta \parallel^{*}$
gives the support for $Nf$, where the norm is such that $H'(E) \rightarrow Exp E^{*}$. Consider the type for $\widehat{\widehat{Nf}}=\widehat{\widehat{f}(\psi^{*}-i\eta)}$, where $\eta$ is real. If we have $\sup \mid e^{<\eta,\psi> \widehat{\widehat{f}}}(\psi) \mid \leq c$, then $\sup \mid \eta \mid$ can be used as a norm. Given that $\sup e^{\mid \eta \mid \mid \psi \mid}
\mid \widehat{\widehat{f}} \mid \leq c$ and we see that $\widehat{\widehat{f}}$ has negative order $\mid \eta \mid$.

\section{Fundamental representations}
Naturally, a discussion on the global models is dependent on the representation
of the symbol.

Consider $\big[ e^{\Phi + C}f \big]^{\widehat{}}=\big[ e^{< \eta,\psi>}f \big]^{\widehat{}}$, where
$\Phi=\tilde{L}(\psi)$ for a homomorphism $L(e^{\psi})=e^{\tilde{L}(\psi)}$. We assume that
 $<\eta,\psi>$ is a scalar product in the sense that
 $< \alpha \eta,\beta \psi>=\alpha \overline{\beta}< \eta,\psi>$, that is $<\eta,-\psi>=-<\eta,\psi>$,
why the phase in the right hand side is odd and we assume that $\tilde{L}$ is odd. We have earlier noted (\cite{Dahn13})
presence of unbounded sublevel sets $\{ \phi < C \}$ or better the condition 
$< \eta,\psi><0$ in the right hand side. will give raise to continuum or functional representations
in the spectral resolution. 
 Note that the condition on $f_{1} \bot f_{2}$ in the real space, 
does not necessarily imply the same orthogonality for $e^{< \eta,\psi>}f$. According to Wiener's theorem, any 
function $g \in L^{1}$, can be approximated by translates of $f$ iff $\{ \tau \widehat{f} \neq 0 \}$ and $f(0) \neq 0$. Thus, if $\widehat{f}  \neq 0$ and $(e^{<\eta,\psi>}f)^{\widehat{}} \neq 0$, we can prove that the orthogonality 
is preserved. 

\vsp

Assume existence of $\eta$ such that $<\eta,g>=0$ for $g \in E_{0} \subset C$ implies $<\eta,g>=0$ for $g \in C$, so that $E_{0}$ is approximating $C$ uniformly. For instance $f=e^{\phi}$ if $\{ \phi - c \}=\{ < \eta,g>=0 \}$, for a constant $c$. If the moment problem is discussed for the phase, such that we are given $\{ \mu_{n} \}$ real constants and $\mu_{n}=<\eta_{n},g>$, for some $\eta_{n}$, that is $\mu_{n}=\mu_{n}I$ and $\eta_{n}$ gives a resolution of the identity. 

More generally, we can give an approximation property by existence of $L F_{T}=\gamma_{T}$ and $E_{0}$ such that $LE_{0}=0$ and some $f_{0} \notin E_{0}$ such that $Lf_{0}=1$. In this case $E_{0}$ approximates the continuous functions. We can define $E_{0}$ as the set where the phase $\phi$ has representation by a scalar product $<>$ for some kernel $\eta$, that is $\phi \neq <>$ iff $\exists L^{-1}$. In this case $L^{-1}$ maps $<> \rightarrow \phi$, that is the phase can be reconstructed by the approximation property. Note finally, that if the kernel representing the operator is symmetric in $x,y$, it is sufficient to consider real valued phases, for approximations of complex valued phases in a finitely generated symmetric ideal.

\subsection*{The representation theorem}
Consider the following problem, when for a continuous homomorphism $L$ and 
$\mathcal{L}'=\{ \exists \mbox{!  } \eta \quad L(e^{\psi})=e^{< \eta,\psi>} \}$, 
do we have $L \in \mathcal{L}'$. Let $\mathcal{L}^{0}=\{ \exists \mbox{!  } 
\eta \quad L(\psi)=< \eta,\psi > \}$, where in this case $L$ is assumed 
continuous and linear. If $X=H(\Omega)$, for an open set $\Omega$, 
we assume $\widehat{X}=\{ e^{\psi} \quad \psi \in H(\Omega) \}$, then 
$L \in \mathcal{L}^{0}(\widehat{X})$. If we have existence of $\eta_{x}$, 
for $x$ fix, we have $L \in (\widehat{X})'$. Assume for $M \in (X)'$, that 
$< M,\psi>=<L,\widehat{\psi}>$, then $L=\mathcal{F}^{-1}M$. When 
$L \in \mathcal{L}'$ and if $L$ is algebraic in $e^{\psi}$, that is linear 
in $\psi$, we have $L(e^{\psi_{1}+\psi_{2}})=L(e^{\psi_{1}})L(e^{\psi_{2}})=
e^{<\eta,\psi_{1}>}e^{<\eta,\psi_{2}>}$. Further if 
$L_{1},L_{2} \in \mathcal{L}'$, we have 
$L_{1}L_{2}(e^{\psi})=e^{<\eta_{1}+\eta_{2},\psi>}=
e^{<\eta_{1},\psi>}e^{<\eta_{2},\psi>}$. If 
$\big[ \widehat{I},N \big](\psi)=e^{<\eta,\psi>}$ and 
$\big[ \widehat{N},I \big](\psi)=<\eta,\psi>$. Assume the commutator 
$C$ such that $C \big[ \widehat{I},N \big]=\big[ \widehat{N},I \big]$, 
then $\mathcal{F}^{-1} C \mathcal{F} \big[ I,N \big]=\big[ N,I \big]$.

\vsp

Note that $\big[ \widehat{\widehat{I}},N \big](\psi)=e^{e^{<\eta,\psi>}}$ 
and $\big[ \widehat{I},\widehat{N} \big](\psi)=e^{N(e^{\psi})}=
e^{<\eta,e^{\psi}>}$ $(=e^{<e^{\eta},\psi>})$ and 
$\big[ I, \widehat{\widehat{N}} \big]=N(e^{e^{\psi}})=<\eta,e^{e^{\psi}}>$.
If $\big[ \widehat{I},N \big]=\big[ \widehat{N},I \big]$, we say that $N$ 
is algebraic. Let $<N(\psi),\theta>=<\psi,{}^{t}N(\theta)>$. If 
$N(\psi) \in (\widehat{X})'$ implies ${}^{t}N(e^{\theta}) \in X'$, for 
$x$ fix. Then ${}^{t}N \widehat{I}=\big[ \mathcal{F} N \big]^{t}$. If 
$<\widehat{I}(\psi),\theta>=<\psi,I {}^{t} \mathcal{F} \theta >$ iff 
$< e^{\psi},\theta>=<\psi,{}^{t} \mathcal{F} \theta>$. Let 
$C \big[ N,\widehat{I} \big]=\big[ \widehat{I},N \big]$. If 
$<\eta_{x},\psi> \in X$ implies $e^{< \eta_{x},\psi>} \in \widehat{X}$.
If we assume $N(e^{\psi})>0$, that is $\widehat{N}(\psi)>0$ implies 
$\theta \in X$, such that $N(e^{\psi})=e^{\theta} \in \widehat{X}$.

\subsubsection*{Continuation representations}

The continuations will be related to the representation theorem. Assume for this reason 
$Nf$ locally injective and that $\frac{<\eta_{T},x>}{\phi(x)}=const$ means $T=0$,
for $\eta_{T} \in \mathcal{L}$, that is $\phi(x) \neq < \eta_{T},x>$ locally.
Alternatively $<\eta_{T},x> \leq \phi(x) \leq \leq <\lambda_{T},x>$, where $\eta_{T},\lambda_{T} \in \mathcal{L}$. Over the lineality we define continuations $(A)$ on the form 
$\{ e^{\phi} \equiv 1 \}$ where $\eta \in (J)^{\bot}$ and $(J)$ is the 
ideal of functions with phase $\eta \geq 0$ and $nbhd(\Delta)=\{ \eta \quad <\eta,\psi> \geq 0 \}$.
If $c_{1} \leq \frac{A P}{P} \leq c_{2}$ in $\infty$ then $AP \sim P$.
If $P$ is hypoelliptic with respect to $\Delta$, we have for every $A \in (A)$ that $A P \neq P$.
Note the interpolation problem to determine $v_{1}$ so that given the arithmetic mean $M$ is continuous, we have
$e^{M(\phi) - \phi}=e^{v_{1}}$, further existence of $v_{2}$ and so on. 
Note also under the condition that $f \bot g$, in the sense that the quotient $f/g \rightarrow 0$ in the 
infinity, if the quotient is algebraic, the lineality set containing the infinity,
is standard complexified.

\vsp

In the case with non finite Dirichlet integral, we consider mainly three cases.
Assume $f=e^{\phi} \notin L^{1}$ but
$\phi \in L^{1}$, or $e^{\phi + v} \in L^{1}$ $\exists v$, or $\phi/<> \in L^{1}$
In the last case we note that if $\phi/<v,z> - <w,z> \rightarrow 0$ and if $<v,z> \neq 0$
in a neighborhood of the infinity, we have that $\phi \sim <v,z><w,z>$ near $\infty$.
If we consider the inverse mapping
$z(v)=<v,z>$ we have that $z(v)z(w) \sim z(vw) \sim z(\phi)$ can be seen as a condition on algebraicity
for the inverse. Note that the condition $\phi - <> \rightarrow 0$ in $\infty$ means that
$\phi \sim <> + P(\frac{1}{z})$ in $\infty$, assuming the difference analytic.
A sufficient condition should be that $1/e^{\phi-<>}$ is reduced in $\infty$..

\section{The lifting principle and the moment problem}
The problem of finding a holomorphic $F_{T}$ such that over trajectories to a dynamical system $\gamma$,
$F_{T}(\gamma)(\zeta)=f_{T}(\zeta)$, can be compared with the problem of finding a Stieltjes measure
$v_{T}$ solving the moment problem, under suitable conditions on the ramifier.  We consider $\varphi(f)(\zeta_{T})=\zeta_{T}$ and $\vartheta(\log f)(\zeta_{T}) =\zeta_{T}$, when $\log f \in L^{1}$.

The problem of determining the spiral cluster sets is in most cases unsolved
when the Dirichlet integral is not finite.
Denote, $H(\alpha)=\int \frac{d \phi(t)}{\alpha - t}$, for $\mbox{ Im }\alpha \neq 0$
Then the moment problem is determined iff $H(\alpha)$ is constant. Over solutions $\phi$ to the moment problem, we have that $\{ H(\alpha) \}$ is convex (\cite{Riesz23}).
We can write $\int F(t) d V_{\lambda}(t) \sim (V_{\lambda} F)=\int K_{\lambda} F d t$, where $K_{\lambda}$ is of Carleman type, then for these kernels $0$ is in the limit for spectrum (\cite{Riesz56})

\vsp

 Let $\Delta$ be defined through $Tf(y+i t \eta)=e^{<\eta>} Tf=Tf$, where $f$ is analytic or 
equivalently $<\eta> \equiv 0$, relative the transform $T$. Let $E=\Gamma^{\bot}$. The moment problem can now be compared with the possibility 
of uniform approximation by $\Gamma$ through $\Gamma^{\bot}$ or equivalently the proposition that there does not 
exist $d \alpha$ of bounded variation such that $Tg \bot d \alpha$ $\forall Tg \in E$ or unique 
existence of a functional $N$, such that $NTg \equiv 0$ (possibly modulo $C^{\infty}$) on $E$ and $NTf \equiv 1$ (modulo $C^{\infty}$) on $\Gamma$. In this context 
a determined moment problem means presence of a normal model.

Assume $f=f_{1} + i f_{2}$ such that $T f_{1}=e^{\phi}T f_{2}$, where $0 < \phi \rightarrow \infty$ on $\Gamma^{\bot}$. Further that $e^{\phi}T=T \sigma$ in $L_{loc}^{1}$.
Assume $M$ defines invariance in $L_{loc}^{1}$, when $\phi$ a constant. If the moment problem is determined as above, we may have that $M \equiv N$ on $\Gamma^{\bot}$ and $MTf=1$. Given a determined moment problem there can thus be an intersection of sets of invariance.

\subsection*{Interpolation property, approximation property}

Assume $\frac{N_{1}f}{f}=e^{<\eta_{1},\phi>}$. Then an interpolation property is existence of $v_{1}$ such that
$e^{M(\phi) - \phi}=e^{v_{1}}$, where $M$ is the arithmetic mean. If we let $v_{2}=M_{2}(\phi)-M(\phi)+M(\phi)-\phi$, then $M(v_{1})+v_{1}=v_{2}$. In the same manner for $v_{3}$. If we assume $M$ continuous and single valued, existence of $v_{1}$ implies existence of $v_{2}$ and so on. If we assume $M_{N}$ continuous for $N \geq N_{0}$, then $M_{N}(\phi)-\phi=v_{N}$ does not imply $M_{N}(\phi)-M_{N-1}(\phi)+\ldots=v_{N}$.
The main problem is to determine how the interpolation property relates to presence of clustersets relating to lineality and orthogonality. Also to determine how the clustersets depend on the lineality and orthogonality. Let $N_{1}f=e^{\phi + v_{1}}$, that is $f$ can be continued to $e^{M(\phi)}$, If $N_{1}'f=e^{\phi + v_{1}'}$, then we have existence of $N_{1}''$ such that $\phi + v_{1}' + v_{1}''=\phi + v_{1}$.

Assume $f, \mbox{ Re }f \in (I) \subset L^{1}$ $(L^{q})$, where $(I)$ is assumed without lineality. We then have existence of an integer $N$ such that $(\mbox{ Im }f)^{N} \in (I)$. If $f=e^{\phi}$ and $\overline{f}=e^{\overline{\phi}}$ with $\overline{\phi}(z)=\phi(\overline{z})$ (locally)
If $e^{\Sigma \phi_{j}} \in (I)$ (algebraic), then $e^{\phi_{1} + v} \in (I)$, where $v=\Sigma_{2}^{N} \phi_{j}$. Thus orthogonality corresponds to an interpolation property over $(I)$. 

For instance if $f \in (I) \subset L^{1}$ and $g \bot f$, then we have existence of $v$ such that 
$v g \in (I)$. If $f \mid_{V} \equiv 0$,$V=N(I)$ and $g \bot f$ on $f \neq 0$, then there is a $v (\sim f)$ such that $vg \mid_{V} \equiv 0$..

Assume $\Omega_{\mu}=\{ y \quad T f(y) < \mu \}$, then $y + \eta \in \Omega_{\mu}$, for $\eta \in L$ a 
line containing $\infty$. 
In presence of lineality, we are looking for existence of $\eta$ such that $< \eta, \psi > \equiv 0$, in this case the sets
$\Omega_{\mu}$ are unbounded.
When $\Omega_{\mu} \subset \subset \Omega$, then $\int_{\Omega_{\mu}} e^{<\eta,\cdot>} d v$ defines a regular function, under suitable conditions on $\eta$.

\vsp

Note that if $1/\rho$ is finite, we have that $\{ \mu < \rho \}$ is unbounded, for small $\mu$
Assume $1/\rho(z,\overline{z})$ entire, this can be related to a condition on finite order.
For instance if for $f$ entire, $f(\frac{1}{z})=\rho(z,\frac{1}{z}) \frac{1}{f(z)}$, then $f$
is not reduced in both $0$ and infinity simultaneously. Note that if $S_{j}$ is the set where $\rho(z,\overline{z})=c_{j}$
a constant, then we have convergence $S_{j} \rightarrow S_{0}$, where $c_{0}$ is the constant defined by 
$c_{j} \rightarrow c_{0}$ (\cite{Nishino68}).

\subsubsection*{The separation property}

We should also discuss how the separation problem is related to the interpolation property.
Assume $A,B$ the disjoint supports of two functions $a,b$ and that $\int a \overline{b}=0$. 
The interpolation property is that we have existence of a function $c$ with support a line $L$ (oriented) 
that separates $A,B$. For instance we can write $\gamma = \overline{\gamma}$ iff $c(x) \neq 0$. 
That is, assume existence of $L$ between $A,B$ and that $L$ is oriented. Given a point (segment) on $bdA$ we 
can determine a line $l$ (oriented) between this point and a point (segment) on $bdB$ and that passes $L$ 
from one side (-) to the other (+). Let $\mathcal{U}$ denote the space of all lines $l$, formed in this way 
such that $a \overline{b}=0$ on $\mathcal{U}$. If we assume $a,b$ analytic functions, we can assume 
$\mathcal{U}$ a domain of holomorphy. On any line $l$ between $A$ and $B$, there is a point where 
$c(x) \neq 0$. We can now consider the supports as one sided (possibly in different leaves.). 
An example is the condition $a \prec b$, if $a/b \rightarrow 0$ in the infinity, and existence of $c$ such that 
$a \prec c \prec b$. If for instance $\mid c \mid=1$, we have $\inf a \overline{c} c \overline{b}=0$, 
further if $m \Omega > 0$ and $L \cap \gamma = \emptyset$, then $\int_{\gamma} a \overline{c}=0$.
Note also that the proposition that $a/c$ is algebraic locally, is implied by a finite distance to essential
singularities, given that $a,c$ are entire.

\subsubsection*{Monogenity and the interpolation property}
The condition on monogenity is intended to give the biggest possible class of trajectories for which 
continuation is possible. If $V$ is the domain where $v \in L^{1}$, we can denote $V^{*}=\log V$.
That is if $v \in L^{1}(V)$, then $\widehat{v} \rightarrow 0$ on $V^{*}$. We can assume $V$ 
unbounded, that is $e^{\psi} \rightarrow \infty$, when $\widehat{v}(\psi) \rightarrow 0$. If we assume 
$\widehat{I} \widehat{\eta}(\psi)=\widehat{\eta}\widehat{I}(\psi)$, then we are assuming $\psi \in V^{**}$, 
that is $e^{e^{\psi}} \in V$. 
\vsp

Assume $e^{<\eta,\psi>}f=v(e^{e^{\psi}})$, where $v \in L^{1}$ (or $v \in \mathcal{D}_{L^{1}}'$), 
that is we are assuming that $e^{e^{\psi}}$ is in the domain for $v \in L^{1}$, which is perceived as the outer 
frame for continuation. If $\widehat{N}(\psi)=e^{<\eta,\psi>}=\widehat{I}\widehat{v}$, where 
$\widehat{v}(\psi)=<\eta,\psi>$. Thus, if $\widehat{I}\widehat{v}=\widehat{v}\widehat{I}$, we have 
$e^{<\eta,\psi>}=v(e^{e^{\psi}})$, for $v \in L^{1}$. Consider now for the means $M_{N}$, 
$A_{1}=\mathcal{F} M_{1} \mathcal{F}^{-1}$, $A_{2}=A_{1}^{2}$ and so on, that is 
$M \mathcal{F}^{-1} \eta_{1} \sim v \in L^{1}$. Let $\eta_{N}=\mathcal{F} M_{-N} v$, we then have that 
$A_{N} \eta_{N} \sim \mathcal{F} M_{N} \mathcal{F}^{-1} \eta_{N} \sim \eta_{0}$. If we assume 
$\mathcal{F}^{-1} \eta_{N}$ has support in $\{ 0 \}$ and as $\widehat{v} \rightarrow 0$, in the infinity, 
that is $\mid \widehat{v} \mid \leq C/\mid T \mid^{N}$, as $\mid T \mid \rightarrow \infty$ and we have 
$\mid \eta_{N} \mid \leq \epsilon$, as $\mid T \mid \rightarrow \infty$. Note particularly, that 
$Nf-f=0$ iff $<\eta,\psi>=0$. Finally note the example, where $\mid T^{N} \eta \mid \leq C$, as 
$\mid T \mid \rightarrow \infty$, that is $\mid \widehat{(\frac{d^{N}}{d T^{N}}v)} \mid \leq C$, as 
$\mid T \mid \rightarrow \infty$, that is $\frac{d^{N}}{d T^{N}}v$ has support in $\{ 0 \}$.
Assume now $\Phi=\tilde{L}(\psi)$ and $\Phi + C_{0}=<\eta_{T},\psi>$ and $\Phi_{1}=M(\Phi)$ ($M$ is the arithmetic mean) with 
$\Phi_{1} + C_{1}=<\eta_{T}^{1},\psi>$. Then for some $N$, we have $\Phi_{N} + C_{N}=0$ and
that $\mbox{ supp } \eta_{T}^{N}=\{ 0 \}$. Note that when $L$ is algebraic it maps Puiseux series onto Puiseux series. 

\section{Comparison with the spectral theory}

Assume the spectrum $\sigma$ is the set where we do not have representation using a regular analytic function. The condition $h(x)/x < \mu$ close to the infinity,
for a constant $\mu$, means that the domain $(x,h(x))$ for one-sidedness is unbounded, and the resolution
over this domain is represented by a functional. If we consider the condition $h(x)/x \equiv \lambda$
on a set of positive measure (spectrum) for $\mid x \mid=1$, we have $h(x)/x \equiv \lambda$,
for $\mid x \mid < 1$, that is for $1/x$ large.
Compare also the theory of singular integrals, $H(x)=\frac{1}{x}h(\frac{1}{x})$ on $0 < \mid x \mid < R$
where $H$ is regular. If we consider the mapping $\Gamma_{0}=(X_{0},Y_{0}) \rightarrow \Gamma_{1}$
according $\Gamma_{1}/\Gamma_{0}=\rho$, we consider the condition $\{ \rho < \lambda \}$
bounded, locally in a small neighborhood of a point. If for instance $X_{1}=\rho_{1} X_{0}$,
we have that $X_{1}=\rho ' X_{0} + \rho X_{0}$, thus we have that $\{ \rho ' = \rho = 0 \}$
has isolated points. If we start from the representation $e^{\phi_{1} - \phi_{0}}$, where
$\phi_{1}-\phi_{0} \in L^{1}$. If $\int_{V} e^{\phi_{1} - \phi_{0} - <>}d \sigma=0$,
we have that either the integrand $\equiv 0$ on $V$ or we have $\sigma(V) =0$.
Note that if $T$ is a mapping such that $-T^{2}=I$ we have the condition $\frac{Tf}{f} < \mu$
is a condition on lower boundedness $\frac{g}{-Tg} < \mu$, where $\mu$ is a (negative) constant.
Note that under the condition that the phase is in $L^{1}$, we have finite order singularities,(\cite{Riesz15})
why the complement to this set can always can be chosen as regular. We can prove that given a Tauberian condition on the phase, (\cite{Dahn13}) we can under the condition of summable phase, prove existence of a regular approximation of
any singular point.

\newtheorem{prop1}{Proposition}[section]
\begin{prop1}
 The spectral function corresponding to a self-adjoint realization of a hypoelliptic, constant coefficients differential operator $P(D)$ on $L^{2}(\mathbf{R}^{n})$
can be represented by a regular kernel
$$ e_{\lambda}(x,y)=\frac{1}{(2 \pi)^{n}} \int_{P(\xi) < \lambda} e^{i < x-y,\xi>} d \xi$$
$e_{\lambda} \in C^{\infty}(\mathbf{R}^{n} \times \mathbf{R}^{n})$
 The spectral kernel corresponding to a self-adjoint realization of a constant coefficients differential operator in the radical to ideal of hypoelliptic differential operators, can be represented by a functional kernel $e_{\lambda} \in \mathcal{D}'(\mathbf{R}^{n} \times \mathbf{R}^{n})$.
\end{prop1}
(\cite{Nilsson72},\cite{Dahn13}) The results can be generalized to variable coefficients formally hypoelliptic differential operators.

\subsubsection*{ The momentproblem}
 For instance if $\mu_{n}=(T_{n} x ,y)=\int f_{n} d \mu_{x,y}$ (\cite{Riesz23}), we have existence of ${}^{t} S$ such that if $\mu_{n}=0$ for all $n$, $({}^{t} S T x, y)=0$, for all $x,y$ that is $(Tx,Sy)=0$ or $S \bot T$, for some ${}^{t}S$. The representation $<f_{n}, d \alpha>=\mu_{n}$, corresponds to weak convergence, (\cite{Riesz56}). 
We do not exclude the case when $(d \alpha)^{\bot}$ is non trivial and of infinite order.
When $(Tx,y)=\mu_{n}(x,y)$, then $\mu_{n}$ are spectrum for the operator $T$ and in the momentproblem the growth points. If we assume $(x,y)(\zeta)$ defines a spiral, we are considering $\zeta \rightarrow z \rightarrow \mid z \mid (\zeta)$. Spirals can be defined using implicit functions, as in $u(\eta(\zeta))=\widehat{v}(\zeta)$, where $\eta$ is locally 1-1. Assume $\phi(\zeta,z_{1})$ a polynomial in $z_{1}$, where $u$ is holomorphic on $\{ F_{1}=0 \}$ and $u=\phi/\frac{\delta F_{1}}{\delta z_{1}}$. If $u$ is holomorphic on $\mid \zeta \mid < A$ and $\Sigma$ is $z_{1}=\eta_{i}(\zeta)$, where $\eta_{i}$ is real analytic, we consider $(\zeta,\eta(\zeta))$ a regular path. We assume $F_{1}(\zeta,z_{1})=0$ implies $d_{z_{1}}F_{1} \neq 0$ and $u(z) d_{z} F_{1}=\Sigma c_{j}(\zeta) z_{1}^{j}$ on an irreducible component in $\{ F_{1}=0 \}$ and we have $u(\eta(\zeta)) d F_{1}(\zeta,\eta_{1})=\Sigma b_{j}(\zeta) \eta_{1}^{j/\mu} $.
We also note the following example. Assume $\widehat{\mu}=\int e^{-t x} d \mu$. 
Assume $M_{N}(\widehat{\mu})=c_{N}$ constant, where $M_{N}$ is the iterated arithmetic mean and 
$\widehat{\mu}+\eta_{N}=c_{N}$ with $\eta_{N} \geq 0$ decreasing implies $
\widehat{\mu} \leq c_{N}$, where $\eta_{N}$ is bounded for $N$ large. 
Consequently the set $\{ \eta_{N} > \lambda_{n} \}$, where $\lambda_{n}$ is
a fixed constant is a decreasing and bounded set, why we can use function theory for large $N$. Assume $\widehat{\mu}(x_{0})=c_{0}$ and 
$c_{1}=M(\widehat{\mu})(x_{1})$. Determine a path between $x_{0}$ and 
$x_{1}$. If $c_{0}+\eta_{0}=\widehat{\mu}(x_{0})$, then $\eta_{0}=0$ in
$x_{0}$. In the same manner, if $\widehat{\mu}+\eta_{1}=c_{1}$. Further, if 
$\eta_{2}=\eta_{1} + M(\eta_{1})$, this implies 
$\eta_{1} \geq M(\eta_{1})$.

\subsubsection*{Continuation with disjoint support}

The approach is using $F_{T}=e^{\phi_{T}}$, where $\phi_{T}=<> + \Sigma v_{j}$, according to the 
interpolation property. The problem is if $v_{j}$ can be given with disjoint support, 
respecting the interpolation property. If for the vorticity, $W(\phi_{T} - \Sigma v_{j})=0$ on a set 
$\Omega_{0}$, then the infinitesimal movement can be given as translation (symmetry condition). 
Further, we can consider $\Omega_{+}=\{ W(\phi_{T} - \Sigma v_{j}) > 0)$ and analogously for $\Omega_{-}$.

\vsp

Consider now a disjoint decomposition of spectrum from a group $\mathcal{G}$.
Assume $g=q$ a polynomial such that $\frac{1}{q^{(j)}} \rightarrow 0$, $j=0,\ldots,n-1$ and 
$T \rightarrow \infty$. We assume $f_{1/T} \sim \frac{1}{q_{T}}$. Then
$\frac{1}{q^{n-1}}\frac{q^{n-1}}{q^{n}} \rightarrow 0$, as $T \rightarrow \infty$, since $q$ polynomials. 
Thus, if $\frac{1}{q^{j}}=p_{j}$ and $x q^{j-1} \sim q^{j}$, we have 
$\frac{1}{q^{j}} \sim x^{j} p_{j} \rightarrow 0$, as $T \rightarrow \infty$. This means that for 
$j=0,1,\ldots,n$, if $q^{j}=const.$ we have $1/T=0$. In this sense if $q^{j}$ quasi orthogonal, we can chose 
$p_{j}=\widehat{v}_{j}$, for $v_{j} \in L^{1}$. Conversely, assume $D^{\alpha}f_{\alpha}$ with 
$f_{\alpha} \in L^{1}$, that is $\widehat{f_{\alpha}} \rightarrow 0$. Assume now that 
$\xi^{\alpha} \widehat{f_{\alpha}} \rightarrow 0$, as $\mid \alpha \mid \leq m$, for $\widehat{f_{\alpha}}$ 
corresponding to quasi orthogonal functions.Characteristic for quasi orthogonal polynomials, of the same degree is that given that their quotient 
is non constant, they do not have common zero's. If the coefficients are non real, the zero's are in a 
half space (\cite{Riesz23}). Existence of base must depend on singularities.
If the spectrum is the support of the polynomials, then the for quasiorthogonal polynomials, we have 
a disjoint decomposition of the spectrum. This means that the spectrum is algebraic. For example assume 
$e \in G$ and $p,I \in G$ and $pe \equiv I$ with $\sigma(p) \cap \sigma(e)=\{ 0 \}$. This means that we can not 
find $c$ with support on a separating line, such that $e < c < p$. 
Conversely, if we have existence of $c$ as above, then $0 < \mid p - \frac{1}{e} \mid < \infty$. 
If $p,q \in G$, with $pe=qf=I$ and where $pq$ is polynomial and if $\sigma(e) \cap \sigma(f)=\emptyset$, then 
$ef \equiv 0$ on $\Omega$. Thus, we must have $\sigma(e) \cap \sigma(f)=\{ 0 \}$ and $p \bot q$. If $e < c < 1$,
where $c$ has support on a separating line, that is $e$ has one sided support. Let us chose $q$ such that 
$q(x)=0$ iff $e(x) \neq 0$ and $\Omega$ is such that $qe=1$ on $\Omega$. If $p,q \in G$ have onesided supports 
with $\sigma(e) \cap \sigma(f)=\{ 0 \}$, we have an interpolation property and conversely.

Another example is $\int_{\Gamma} pq d \varphi=0$, where the integrand is polynomial (does not mean $p$ 
polynomial), where $\Gamma$ is semialgebraic or semianalytic $\{ p \geq 0,q \geq 0 \}$.

If we assume $[I,h]=[h,I]$ with compatibility conditions, assume $h=h_{1}h_{2}$,
where one of them is even and the other odd, then $h$ is even. If both are odd,
then $h$ is odd. Assume $h$ even, then $dh$ is odd and $g$ odd means $d g$ even.
In this case $d (fg)(-x)=(- gdf + f d g)(x)$, that is $d (fg)(-x)=- d (fg)$
means that $f d g=0$.

\subsubsection*{Determined tangents in phase space}
We have noted that a condition for a ps.d.o representation, is that we have tangent lines $h(f)=\mu f$ (\cite{Dahn13}). We would like a representation, when there are no tangent lines present.
Assume $h(e^{\varphi})=e^{\tilde{h}(\varphi)}$, then the condition $\tilde{h}(\varphi)/\varphi=\mu$, for a constant $\mu$ means that $\frac{d}{d T} \frac{\tilde{h}(\varphi)}{\varphi}=0$. Assuming $\tilde{h}(\varphi)$ linear in $\varphi$, we can then write $\tilde{h}(\varphi)=<\eta.\varphi>$, for some $\eta (\in H')$, where the tangent is given by $\mbox{ supp } \eta=\{ 0 \}$. Assume for simplicity that $\frac{d \tilde{h}(\varphi)}{d \varphi}=\frac{d \tilde{h}(\varphi)/ dT}{d \varphi/d T}$. Let $g_{\mu}(\varphi)=\frac{\tilde{h}(\varphi)}{\varphi}-\mu$, for a constant $\mu$. Consider the ideal $E$ where $g_{\mu}(\varphi)=g_{\mu}(\frac{{}^{t} d \varphi}{d T})=0$ and $(J)=(\mbox{ ker }g_{\mu})$. 
We have $\varphi \in J$ iff ${}^{t} \frac{d}{d T} \varphi \in J$. Thus, if the kernel corresponding to $g_{\mu} + g_{\mu} \frac{d}{d T}$ is $\eta$, then $\mbox{ supp } \eta=\{ 0 \}$. If $J$ is given by $<\eta,\cdot>=0$, then using the nullstellensatz, we have existence of a $\mu$ such that $\varphi^{\mu} \sim < \eta,\varphi>$..
If we assume the ideal $(J)$ is radical, we have $\eta \sim \delta_{0}$. Note that if the two lines $L_{j}(\varphi)=<\eta,\varphi>-\mu_{j} <\delta,\varphi>$, $j=1,2$, then $L_{1}(\varphi)=L_{2}(\varphi)$ iff $\mu_{1}=\mu_{2}$. 

\newtheorem{cont}[prop1]{Remark}
\begin{cont}
 The continuations that we are considering are assumed with representation in $\mathcal{D}_{L^{1}}'$ and with
a tangent determined, that is we assume presence of tangent lines in phase space.
\end{cont}

\subsubsection*{Parseval's relation}
The Parselval relation is given $\mid f \mid,\mid f \mid^{2} \in L^{1}(d \phi)$ we have that the moment problem is determined or $\phi$ is extremal. In this case we have that $f(x)$ has Fourier series representation in growth points for $\phi$.
Consider $Af(x)-\mu f(x)=0$ and $t \rightarrow \lambda(t)=x$. Starting from a continuous curve $S=\{ t=\varphi(x) \}$
with a jump on $\{ x = const \}$. Then we have $\int_{\alpha}^{\beta} F(t) d t=\int_{a}^{b} f(x) d \varphi(x)$,
with $F(t)=f(x(t))$.
Assume $(I)=\{ f \qquad d h(f)=0 \}$ (implies $f=const$) and $(J)=\{ g \notin (I) g$ regular and locally algebraic $\}$. The problem is to determine when $f \in (I) \bigoplus (J)$. Given that $\Sigma \mid a_{n} \mid^{2}=\int \mid f \mid^{2} d \varphi$, we can find a polynomial $s$ of degree $\leq n$, that minimizes $\int \mid f - s \mid^{2} d \varphi$ and in this case the transversal can be given as a locally algebraic function. Assume we have existence of $h_{2} \in C^{1}$, such that $\frac{d h_{1}}{d z}=h_{2}$- In general we have that $\mbox{ supp }h_{2} \subset \mbox{ supp }h_{1}$. If now $h_{2}$ is reduced, we have that we are done.
If we assume $h_{0}$ has compact support, we have that $h_{n} \in \mathcal{E}^{(0)'}$. If $I_{1}=\{ h_{1}(f)=0 \} \subset I_{2} \subset \ldots$ and $f_{0} \in I_{2} \backslash I_{1}$ a transversal, then $I_{1} \cap I_{2}=\{ 0 \}$ and so on. Thus, if $(I_{2})$ is transversal, then $(I_{1})$ is minimally defined. Further, if $I_{N}$ is transversal, then $f \in \cap_{1}^{N} I_{j}$, then $f=0$. If we assume the transversal is a regular approximation,
$\frac{d}{d z}f_{2} \neq 0$ and $h_{2}=\frac{d}{d z}h_{1}$, thus $f \in I_{2}$ and $\frac{d}{d z}f \neq 0$ implies ${}^{t} \frac{d}{d z}f \in I_{1}$. Note that the condition on the symbol ideal $\frac{d h}{d z} \rightarrow 0$, $\mid z \mid \rightarrow \infty$, means that $h_{1} \rightarrow 0$ as $\mid z \mid \rightarrow \infty$, that is $h_{1} \sim_{m} P(1/z)$ as $\mid z \mid \rightarrow \infty$.

\subsection*{Grouptheoretical considerations}

A general problem in this connection, when the global model is given by a group, is if a group is its center, that is $g \in G$ iff ${}^{t}g \in G$
In this case $F(g \gamma)(\zeta)=F(\gamma)({}^{t} g \zeta)$.  Another problem is the following. Assume the group defines analytic continuation and that $gf=f$ on $\Omega_{g}$,
further that $h f=f$ on $\Omega_{h}$, where $g,h$ are in the group, then obviously $g=h$ on $\Omega_{g} \cap \Omega_{h}$.
That is intersection of sets corresponds to invariance f\"or group action. If all the sets are disjoint, then there is no room for invariance. The problem of determining a group
that can (uniquely) determine the global model in question, can for this reason be related to sets
(or lack of) invariance or isotropy.

\vsp

Assume $g=g_{1}g_{2}$ and $\big[ I,g_{i} \big]=\big[ g_{i}, I \big]$, $i=1,2$, that is if $g_{i}$ is generated by a group,
it is in the center of the group. Then these conditions do not necessarily imply the property for $g$.
When $g \in \mathcal{G}$ is algebraic, we have $\sigma(g f)=\sigma(f)$ (spectral mapping theorem).
The condition
on irreductibility means that if $g=g_{1}g_{2}$, then $g_{1}=id$ or $g_{2}=id$. Note the example,
when $g$ is irreducible, and change of order of integration is motivated, $\big[ g_{1}g_{2} \varphi \big](f)=\big[I,g_{2}\varphi \big](f)=\big[ \big[ I,g_{2} \big],\varphi \big](f)$
and $\varphi(g_{1}g_{2} f)=\big[ \varphi, \big[ I,g_{2} \big] \big](f)$ and if $\big[ I,g_{2} \big]=\big[ g_{2},I \big]$
we have $g \big[ I,\varphi \big]=\big[ \varphi,I \big] g$. We have $\big[ \varphi , gI \big](f)(X)=
\varphi(f)(gX)$ and $\big[ gI,\varphi \big](f)(X)=\varphi(g f)(X)$. More generally, $\big[ \varphi, gI \big]=
\big[ {}^{t}g I,\varphi \big]$ and if $g=g_{1}g_{2}$ is irreducible, $\big[ \varphi , \big[ I,g_{2} \big] \big]=
\big[ \big[ {}^{t}g_{2},I \big],\varphi \big]$. 
However if for a polynomial $P$, $g=P(D)I$ and $P=QR=RQ$, where $R,Q$ polynomials, we can conclude the property for $g$.

\section{Global models}
There are a number of essentially different global models, that can all be combined to a
model adapted to a particular problem. We shall shortly discuss the differences between
the models concerning the possibility of analytic continuation. Given a multiply connected
domain $\Omega$ with a simply connected component $\Omega_{0}$ and $\Omega=\Omega_{0} \cup \cup_{j} V_{j}$.
Assume $F$ analytic on $\Omega_{0}$. We assume that $F$ considered as functional in $\Bm'$
can be continued to $\Omega$. For instance if $F$ is given by a uniform determination on
$\Omega_{0}$, but not on $V_{j}$.
Particularly we can assume a normal model with time axis as transversal. For a continuation
along the transversal to give a reconstruction of an analytic ideal, the ideal should be radical
and such that the mapping from the transversal to the ideal is locally bijective.
An analytic continuation over the time assumes that the movement over time is possible.
In most of the applications we do not assume this. The global models that are based on
clusterset theory, are based on unbounded sublevel sets, that is a continuation as functional.
It is sufficient to assume that $F$ is totally bounded over $V_{j}$. Another example is given
by $f \in I(\Omega)$ such that $f=f_{0}f_{1}$ where $f_{0}$ analytic over $\Omega_{0}$ and
$f_{0} \neq 0$ over $V_{j}$. The possibility for analytic continuation over $V_{j}$ depends
thus on the interpolation property and the property of $f_{1}$ (analytic).

\vsp

Another example is given by $f$ analytic on $\Omega_{0}$ and the time axis as transversal
where $V_{j}$ contains parabolic singularities. We denote $\tilde{f}$ analytic continuation
on $\Omega_{0}$ such that the events that
are given by $V_{j}$ have their maximal points at the transversal and $\tilde{f}$ adapted
so as to describe also these. We are not assuming that $\tilde{f}$ is analytic on $V_{j}$.
We can compare with the condition on a very regular boundary (\cite{Dahn13}).

Also consider a representation at the boundary where $\mbox{ ker }F \neq \{ 0 \}$
and existence of $g_{1}$ such that $\mbox{ ker }(F+g_{1}) \subset \mbox{ ker }F$.
Successively, we can assume $F + \Sigma g_{j}$ has a trivial kernel.
A necessary condition for the possibility of analytic continuation over time in this example,
is that the kernel is trivial.

For the global prediction model on the set $V$, it is necessary to have a representation of the complement of the range $\Omega$ of the symbol to be continued. Note that if $\Omega \subset V$, the complement corresponds to $V \backslash \Omega$, which in algebraic geometry corresponds to the sum of two ideals. The algebraic structure associated to the prediction model is assumed to be a ring in this paper. For a discussion about rings and groups in analytic continuation problems, see (\cite{Nilsson64})

\subsubsection*{Compatibility conditions}
The multiplicative Cousin problem, is given a region $D$ and a covering $\{ U_{j} \}$ with $h_{j}$ 
holomorphic and not identically zero on $U_{j}$, with $h_{i}h_{j}^{-1}$ holomorphic and nonzero on 
$U_{i} \cap U_{j} \neq \emptyset$, we can find $h$ holomorphic, such that $h h_{j}^{-1}$ holomorphic 
and nonzero on $D$. Given Oka's property we have existence of a continuous solution iff we have existence 
of a holomorphic solution. Assume $h_{i}$ such that $\frac{1}{g} h_{i}(\frac{1}{g})=\frac{1}{g} \in U_{i}$. 
If $f_{1/T} \in E_{0}$, a subset of continuous functions, in this application non dense, we have that 
$f_{1/T}$ is uniformly approximating $\frac{1}{g}$, if we have existence of a functional $A$, such that 
$A(\frac{1}{g})=1$ and $A(f_{1/T})=0$. Further, we assume existence of a path for continuous continuation, 
such that $A(\tilde{f_{1/T}})=0$ Assume $E_{0}$ such that $f_{1/T} \neq 0$ and $A(f_{1/T} h_{i}(\frac{1}{g}))=0$.
 If we have parabolic singularities, we can solve the problem using monotropy $A \sim_{m} A_{1}$, with $A_{1}$ 
algebraic. Thus, $A_{1}(f_{1/T} h_{i}(\frac{1}{g}))=A_{1}(f_{1/T})A_{1}h_{i}(\frac{1}{g})=0$, when 
$A_{1}(f_{1/T})=0$.  We say that $f_{1/T}$ has analytic support, if $h(f_{1/T})=0$, for an analytic $h$.  
If $f_{1/T} \in (I)=(\mbox{ ker }h)$, we have $f_{1/T}=0$ on $N(I)$ (boundary).

\vsp

If $\frac{1}{g} h_{i}( \frac{1}{g} )=w(g)=\frac{1}{g}$, we are assuming $g w(g)=1$. Assume for this reason 
$w = \widehat{\mu}$, that is $\widehat{(\frac{d}{d g} \mu)}=\widehat{\delta_{0}}$, that is $\mu$ corresponds 
to Heaviside. If $A(f_{1/T} \widehat{\delta_{0}})=A(\widehat{v}_{T} \widehat{\delta_{0}})=0$, 
$A(v_{T} \widehat{*} \delta)=0$ or $A(\widehat{v}_{T})=0$. Note that if $u(g)=g w(g)$, then $w$ is algebraic iff $u$ is algebraic. Note also that a necessary 
condition for a compatible group system as above, is that $id \sim \delta_{0}$, that is that the supports are disjoint.
If $\mathcal{G}$ is a group for analytic continuation, we note that $g f=f$ in the sense that $\mbox{supp }g f=\mbox{ supp }f$
is not a regular continuation and $g \sim$ id, where id has point support. For a continuation of domain in this case,
we would require an id with continuum support. 

\subsubsection*{Cold case}
Assume $V$ a parabolic set with a maximum point $z_{0}$ with respect to $t$. Assume $V_{1}$ a parabolic event with maximum point $z_{1}$ and assume that the two parabolic sets have points in common or (at least) that we have a bijection $\rho$ such that $\rho_{T}(z)=z_{1}$ and so on until $V_{N}$ which is assumed to be a parabolic event in present time, with maximum point $z_{N}$. We assume the events can be described by dynamical systems (or pseudo bases). The mapping $\rho$ corresponds to a transversal, that is we assume a normal model or even better a normal tube.

With these conditions it is not necessary that we actually have dynamical systems, but we require existence of $\rho$. Thus, $V \rightarrow \gamma_{T}^{N}$ in present time (or $V \rightarrow V_{1} \rightarrow \gamma_{T}^{N}$ or $V \rightarrow \gamma_{T}^{1} \rightarrow \gamma_{T}^{N}$) and solving the problem requires a lifting function $F(\gamma_{T}^{N})(\zeta) \in (I)$. 
We note when the lifting function is used, that it is assumed that $V$ is not in the domain for $F(\gamma_{T}^{N})$, that is we assume existence of $V_{1}$ such that $V \rightarrow V_{1} \rightarrow \gamma_{T}^{N}$.
If we instead are given an event in the future $V$ and $V \rightarrow \gamma_{T} \rightarrow \gamma_{T}^{1} \rightarrow \ldots \rightarrow \gamma_{T}^{N}$ with compatibility conditions, it is possible to solve the prediction problem, given that the respective mappings $\gamma_{T} \rightarrow \ldots \rightarrow \gamma_{T}^{N}$ are continuous.

\subsubsection*{The future}
To do prediction we have to reverse the order and we start with a path in the complement of the range, $\gamma_{T} \rightarrow \ldots V_{T}^{N}$ where a parabolic event can be chosen from the present time or history. If we assume the path can be analytically continued for $F(\eta)$ and $\tilde{\eta}$ along $\gamma_{T}^{j}$, we can look for events $V_{j}$ along the path. Particularly, we can reflect through a path that is not possible to continue analytically and change base and then continue the new path.
The problem of finding a reflecting base can be compared with the reflection through the real axes for the phase, when the problem is to find a given event in the future. We then attempt to reach the same event using analytic continuation of $\eta$. More precisely $F(\eta) \rightarrow \eta \rightarrow \gamma \rightarrow \tilde{\gamma} \rightarrow V$ is parallelled by $ F(\eta) \rightarrow \eta \rightarrow \tilde{\eta} \rightarrow V$. We assume the mapping $\eta \rightarrow \gamma$ a contact transform. If changes of base is used and we use the local lifting function (or functional), we have to consider the transmission property for this mapping.

Note that of fundamental importance for the prediction model is the quality of the reflecting base, that is we need to discuss $V \rightarrow \tilde{\eta}$ that has an inverse with good qualities. For instance if we assume $\log \tilde{\eta}$ summable, for instance if $\frac{d}{d t} \tilde{\eta}=(X,Y)$
we are assuming finite order in the right hand side. For algebraic right hand sides, we obviously have a lifting principle. The lifting principle for analytic right hand sides is dependent on a ramifier $\ram$. In the same manner we must show that the event is possible to reach from $\eta$ with analogous conditions. The continuation is assumed as in Lie, that is $\frac{Y_{\eta}}{X_{\eta}}=\frac{\tilde{Y}_{\gamma}}{\tilde{X}_{\gamma}}=\frac{Y_{\gamma}}{X_{\gamma}}$.

We could speculate that the order of singularities is dependent on the right hand sides of the systems.
However, assuming the continuation of systems satisfy $\frac{B}{A}=\frac{Y}{X}$, repetition of the 
procedure can correspond to a desingularization, that is we can use right hand sides of infinite order $(X,Y)$
and still produce a reflection of finite order $(A,B)$, thus the order of right hand sides are in this context
not dependent of the order of the singularities in $\zeta \in V$.

\subsubsection*{Discussion on reflections}
Assume $\tilde{f}$ the result of reflection through a planar axes, such that $C : \mid \tilde{f} \mid=\mid f \mid$. If two mirrors $R_{1},R_{2}$ are used, then every point on $C$ can be reached by
$L$, a reflecting line relative $R_{1},R_{2}$. The bisectris is the normal to $R_{j}$. Obviously $\mid f \mid \in L^{1}$ iff $\mid \tilde{f} \mid \in L^{1}$. Assume now the first reflection point is $0$ and the second is on $C$. If the reflection result is on the opposite side of the real axes relative $f$, we use two mirrors and if the reflection result on the same side as $f$, we use one mirror. If we assume the reflection points are inner points, any point on $C$ can be reached using only one mirror. If we assume the reflection points are on the axes, then to reach any point on $C$, we have to use two mirrors. 

\vsp

Example: Reflection is not in principle length preserving, but if we assume an interpolation property,
that the reflection point can be moved to the boundary $C$, then the reflection can be realized through simple rotations. If we use one reflection point on the real axes, inner to $C$, then the path $z \rightarrow -z$ is excluded. We use for this case two reflection points inner to $C$ and this means that any point on $C$ can be reflected onto any other point on $C$. The segment on the path between the two mirror points are dependent on the transmission property. We will also discuss a third reflection point in a time parameter. For the applications to prediction, we do not expect the symbol to be symmetric with respect to a reflection point in time. However, it may be analytic independently of the reflection. In higher dimensions, we require several (possibly infinitely many)
reflection points.

\vsp

Assume now $\mathcal{G}_{2}$ the group of reflections through reflection points on the real axes and $\mathcal{G}_{3}$ the group of reflections through reflection points on the axes.
For both groups we assume $\log \mid f \mid \in L^{1}$, which means that we are assuming finite order singularities, that is we are assuming a normal model. Assume $g_{2} \in \mathcal{G}_{2}$ corresponds to a fixed mirror $R_{2} \sim (x,n)$, where $x$ is the reflection point and $n$ gives the angle. Obviously, for every $g_{2} \in \mathcal{G}_{2}$, there is a $g_{3} \in \mathcal{G}_{3}$ such that $g_{2}f=g_{3}f$ (infinitely many).  Further $g f= \zeta f$, for some $\zeta$ with $\mid \zeta \mid =1$. For the group $\mathcal{G}_{2}$, any $f \neq \pm 1$ on $C$ can obviously be reached using $\mathcal{G}_{2}$ in infinitely many ways. However we must except the case when both $f,\tilde{f}$ are on the same axes. These $f$ can be regarded as lacunary for $\mathcal{G}_{2}$. For instance if we assume $x>0$ and $n \leq 0$ and $\mbox{ Im } f \neq 0$ we can reflect onto the lower half plane. For $\mathcal{G}_{3}$, we can repeat the argument. Assume the reflection with two reflection points, then we can map onto any point on $C$ without exception, using two mirrors $R_{1},R_{2}$. Alternatively, assuming parabolic singularities, we can use monotropy, that is the condition that we have existence of a reflection point close to a real point, where normals are $n \geq 0$ or $n \leq 0$ and in this way we can reach all points.

\newtheorem{prop3}[prop1]{Remark}
\begin{prop3}
 Given a point on the boundary to a disk in the plane, we can reach any other point using
reflection and two reflection points on the axes.
\end{prop3}

\subsection*{The reflection model}

Assume $V$ the event we are trying to find information about and $dW$ the base that we use for reflection points 
and tangents. Typically $dW$ is given by a dynamical system or a pseudobase. We will in what follows improperly use the notation $F(X,Y)$ for $F(\gamma)$, that is where $\gamma$ is assumed to satisfy a dynamical system $(X,Y)$.

If we consider $d W_{1}$ at a point in time $A$ and $d W_{2}$ at a point in time $B$, we can consider the 
reflections as overlayers.
At one given point in time, we can in this case give $d W_{1} \rightarrow d W_{2}$ using a contact transform.
For this mapping to be analytic we will in this article require that the events are parabolic. Assume the mapping is $\rho$ and $p$ the mapping between the corresponding events
$V_{1} \rightarrow V_{2}$. When $W_{1} =const$, $V_{1}$ corresponds to foliation. If $Y/X = const.$, we are assuming a tangent determined for the model.

Assume $V_{0} \rightarrow V_{1} \rightarrow \ldots$ are parabolic events that we want to link.
using a bijection $\rho(z_{0})=z_{1}$ that maps max-points
onto max-points. Thus, $(I) \rightarrow W \rightarrow_{p} W_{1}$, $W \rightarrow V$,
$W_{1} \rightarrow V_{1}$ gives $\rho : V \rightarrow V_{1}$, through commutativity.
An example is given by $\gamma \rightarrow \eta$,$\gamma \rightarrow V$, then we can define $V \rightarrow \eta$,
but the mapping cannot necessarily be reversed.

\vsp

Note that the condition that the continuation is pure and planar, in the sense that it has vanishing flux, means that
analyticity is preserved.  However not necessarily continuum, unless the continuation corresponds to a bijection.
When we consider the mapping $V_{1} \rightarrow V_{2}$ and so on in the mirror model,
we note that an absolute continuous mapping preserves zerosets, that is a maximumprincipe
if it is present.

 If $(X,Y) \rightarrow (M,W) \rightarrow (\widehat{M},\widehat{W})=(H,G)$ (\cite{Dahn13}). Let $\Top$ be the completion
of the Fourier-Borel transform, to $-H/G = \Top W / \overline{\Top W}$. Let $F^{\diamondsuit}(H,G)=F(\Top W, \overline{\Top W})$. The condition that $F$ is analytic in its arguments means that $\int_{\Gamma} d F=0$ for all closed contours $\Gamma \sim 0$.

\subsubsection*{Further remarks on the reflection model}
Assume $F(\gamma)(\zeta) \in (I)$ and $G(\eta)(\zeta) \in (I)$, where $G^{-1}$ is defined on $(I)$,
then we have $G^{-1} F(\gamma) \sim \eta$. We are going to discuss 
$F(\gamma) \rightarrow_{G^{-1}} \eta \rightarrow V$, where $\Omega \neq V$. We assume $\gamma$ 
can be continued to $\tilde{\gamma}$, such that $F(\tilde{\gamma}) \in J(U)$, the ideal formed over the set $U$ holding reflecting points and in the same manner 
$G(\tilde{\eta}) \in J(U)$, where we have existence of $G^{-1}$ on $J(U)$ such that 
$\tilde{\eta} \sim G^{-1} F(\tilde{\gamma})$. 

Alternatively, if $f(\zeta) \in (I)$ with base $\{ F_{j} \}$, we are looking for a mapping 
$f(\zeta) \rightarrow \{ W_{j} \}$, where this set can be continued to $J(V)$, the ideal formed over the set $V$. More precisely, we are 
looking for the mapping $W_{j} \rightarrow V$. Thus if $\Omega=N(I)$ with normal $n_{\Omega}$ and $V=N(J)$ 
with normal $n_{V}$, then we assume $n_{\Omega} \cap n_{V} \in bd U$, where we assume $U$ one sided with 
respect to $bd U$. A different problem is to describe $V$, when $U$ is given. In the case when $dW$ is used as a 
mirror, we are exchanging points used as reflection points, to tangents in the point to $W$. Note that $J(U)$ 
is not completely determined by its tangents (blow-up). 
Assume $\Omega$ is one-sided with respect to $bd \Omega$ and define $\Omega^{\bot}=\{ \mbox{ normals } \rightarrow \infty \mbox{ to } \Omega \}$. 
Further, assume $V$ onesided, then we can give $U$ as $\Omega^{\bot} \cap V^{\bot}$. Note that 
$\Omega^{\bot} \cap V^{\bot} \subset \Omega^{c} \cap V^{c}=(\Omega \cup V)^{c}$. Naturally, regularity for this 
set is significant for the definition of $W$. Note that $(I)^{\bot}$ can be defined as $\{ \eta \quad <\eta,\phi>=0 \quad \phi \in (I) \}$, 
that is the annihilators to $(I)$. If we continue to $<\eta , \tilde{\phi}>=0$ and in the same manner 
$<v,\tilde{\psi}>=0$ from $(J)^{\bot}$ and in this way $dW$ can be determined from functionals.
Note that for $\widehat{F}=f$ and $F \in \mathcal{D}_{L^{1}}'$, we can consider the subclass
of $g$ such that $f/g \rightarrow 0$ in the infinity. Using a very regular boundary,
$g$ corresponds to $\tilde{f}_{0}$ in the representation of $f \sim P(D) \tilde{f}_{0}$.
In this sense the boundary representations are the orthogonal to $(I)$.

If $\eta$ is in a dynamical system, we are assuming $F(\gamma) \rightarrow \eta$ maps singularities in the 
foliation to $F$ on singular points for the dynamical system (except for a closed algebraic set). Note that if 
we can chose $G$ algebraic then the same holds for the inverse $G^{-1}$.  Assume $F(\gamma) \in (I)(\Omega)$, 
where $\gamma$ is in a dynamical system with right hand sides $(P,Q)$ and continue $\gamma$ to $\tilde{\gamma}$ 
in a dynamical system with $(X,Y)$ analytical, that defines $W$ $(dW)$.
The tangents $Y/X=const$ corresponds to mirrors or spectrum. For instance $(\mbox{ Re }f, \mbox{ Im }f) \rightarrow (X,Y) \rightarrow \zeta$
and $(\mbox{ Re }f, \mbox{ Im }f) \rightarrow (M,W) \rightarrow \zeta$. If $Y/X$ is analytic, then $\zeta$ 
traces a line. 

Assume $\Omega$ is defined by $\{ f \leq \lambda \}$ and $bd \Omega=\{ f=\lambda \}$, then if $F(\gamma)$ is 
reduced, we have that $\Omega$ is bounded. Note also that for $\tilde{\gamma}$ sufficiently large, we can assume 
existence of $G^{-1}$. If $dW$ is analytic (or meromorphic) on $V$, since $\mathcal{B}_{1} \subset \mathcal{B}_{2}$ implies $\mathcal{B}_{2}' \subset \mathcal{B}_{2}'$, within summable distributions, existence of $G^{-1}$ is not necessarily dependent on involution. However, the transmission property will still depend on the involution condition. 

We conclude with a scheme
\begin{displaymath}
\begin{array}{llll}
I \rightarrow &  I_{1}  &  \rightarrow \ldots &  I_{N}  \\ 
\downarrow & \downarrow &   & \downarrow \\ 
\gamma \rightarrow & \gamma_{1}  &  \rightarrow \ldots & \gamma_{N}  \\ 
\downarrow & \downarrow &   & \downarrow \\ 
\Omega \rightarrow & V_{1}   & \rightarrow \ldots & V 
\end{array}
\end{displaymath}
The prediction model uses $I \rightarrow \gamma \rightarrow \gamma_{1} \rightarrow \ldots \rightarrow \gamma_{N} \rightarrow V$. The mappings $\Omega \rightarrow V_{1} \rightarrow \ldots \rightarrow V$
if they exist can be constructed as a desingularization. Blow-up occurs in the upper row, when the ideals are constructed from the middle row.
\subsection*{The inverse mapping}

 In the reflection model, we assume $T_{1}(f(z))=f(\alpha z)$ corresponds to monotonous continuation from the lower half plane and $(T_{1})_{-1}(f(z))=f(\alpha^{*-1}z)$, continuation to the lower half plane. In the same manner, we define $T_{2}(f(z))=f(z \gamma)$ that corresponds to continuation to the upper half plane. Thus, $\big[ I, T_{1} \big] \sim \big[ T_{2} , I \big]$ assuming both sides are well defined. If we assume $f$ analytic close to the origin, we can find a disk where $f$ is odd. From this disk, we continue the domain to the upper, respective the lower half plane. If the continuation is standard complexified, $z=x + i x^{*}$, where $x$ real. Alternatively, we can consider $T_{1}$ on the phase,
$T_{1} (\log f) \sim T_{1} (\log \overline{f})$, where $\log \mid f \mid \in L^{1}$. An algebraic continuation is standard complexified, without lacunary points, given that the symbol $\rightarrow 0$ on a complex line. If the continuation is not algebraic, we assume the imaginary axes lacunary for the continuation. Examples of continuations that are onesided, are for instance $\tilde{f}$ downward bounded and $\tilde{f}$ is not analytic.
 We consider the the properties $\alpha(-z)=-\alpha(z)$ and $\alpha(-\overline{z})=\overline{\alpha(z)}$.

Compare with essential selfadjointness for the symbol,$F^{*}=\overline{F}$. If $A : T_{1} \rightarrow T_{2}$, with $A T_{1} F = - \overline{T_{1} F}$. Then we have $\mbox{ Re }T_{1} F \sim \frac{1}{2} (I + A) T_{1}F$ and $\mbox{ Im }T_{1}F \sim \frac{1}{2} (I - A) T_{1}F$. If $A^{2}=I$, then $(I + A) \bot (I - A)$. Note that for higher order of transformations, we have $T_{1}f(\gamma_{1})=T_{1}f(\gamma_{2})$, for $\gamma_{1} \neq \gamma_{2}$, that is the inverse mapping $T_{1}f \rightarrow \gamma$ is not locally injective.

We consider the following generalizations of the relation $\big[ I , T_{1} \big] \sim \big[ T_{2} , I \big]$. We consider the symbols where we have that $\big[ I, P T_{1} \big] \sim \big[ Q T_{2} , I \big]$, where $P,Q$ are polynomials. Further the symbols, where we have that $\big[ I, e^{v} T_{1} \big] \sim \big[ e^{\phi} T_{2} , I \big]$, where $\phi,v$ are locally analytic. For example $If(z)=f(z^{*-1})$ and
$\big[ I , P T_{1} \big](f)=I  P(z) T_{1}( \alpha f)=P(z^{*-1})f(z^{*-1} \alpha^{*-1})$,$\big[ Q T_{2}, I \big](f)=Q T_{2}f(z^{*-1})=Q(z) f(z^{*-1} \gamma)$, when $\gamma=\alpha^{*-1}$, which implies $P(1/z^{*}) \sim Q(z)$. If further $T_{1}I=T_{2}I$, we have $T_{1}IT_{2}^{-1}=I$, and we have $T_{1}T_{2}^{-1} \sim Q/P$, that is $T_{1}T_{2}^{-1}f(z)=f(\alpha z \gamma^{-1}) \sim \frac{Q(z)}{P(z)} f(z)$. Further $\big[ I , e^{\phi} T_{1} \big] = \big[ e^{v} T_{2} , I \big]$ and if $\gamma=\alpha^{*-1}$, we have $\phi(z^{*-1}) \sim v(z)$ and
$f(\alpha z \gamma^{-1}) \sim e^{v-\phi}f(z) \sim e^{v(\frac{1}{z^{*}}) - \phi(z)} f(z)$.

\vsp

To determine the inverse $\Phi(f)(\zeta) \rightarrow \zeta$ using Abel's equation (\cite{Julia24}), we must assume isolated singularities. Assume for this reason $\gamma_{T}(\zeta)$ has isolated singularities and assume we allow continuum for $F(\gamma_{T})$ in the first surfaces. Thus, we assume $\gamma_{T}(\zeta)$ has isolated singularities when $\zeta_{T}$ is restricted to a line. With these conditions Abel's equation gives the inverse $\gamma_{T} \rightarrow \zeta_{T}$ and the inverse 
$F(\gamma) \rightarrow \gamma$ is considered to be dependent on involution.

\section{Localization}
Consider the problem of localizing the symbol to a given event $V$ using a mirror (or two). That is given a geometric set $V$ we look for $\gamma_{1}$ (and $\gamma_{2}$) such that $f \in (I)$, can be continued to $\gamma_{1}$ (and $\gamma_{2}$). More precisely if $f \in (I)$ and $f(\zeta)=F(\gamma)(\zeta)$, we assume the continuation is $\tilde{\gamma} : \gamma \rightarrow \gamma_{1} \rightarrow \gamma_{2}$, such that $V$ is in the domain for $F(\tilde{\gamma})$.
 and such that $\tilde{f} \rightarrow V$ continuous. We are assuming a normal model for $\gamma \rightarrow \gamma_{1} \rightarrow \gamma_{2}$. We can further assume $(\gamma_{1})_{T},(\gamma_{2})_{T}$ are dependent on $T$, that is assumed to be dependent on time..
We can take the mirrors at different times $T=t_{1}$ and $T=t_{2}$ and in this case we assume $\mid t_{1} - t_{2} \mid < \delta$, for a small $\delta$. We are assuming $\log f \in L^{1}$ and that $\log \tilde{f}_{T} \in L^{1}$, for all $T$. The continuation is such that $\gamma$ corresponds to a dynamical system with $(X,Y) \rightarrow (X_{1},Y_{1}) \rightarrow (X_{2},Y_{2})$ such that
$\frac{Y}{X}=\frac{Y_{1}}{X_{1}}=\frac{Y_{2}}{X_{2}}=\rho$ (Lie point transformations). Further, that $X_{2}=P$ and $Y_{2}=Q$ polynomials. Thus, $P \rho=Q$, which means that $P \rho$ is polynomial and $\rho$ has finite order, that is the corresponding sections are removable.

The vorticity $W=Y_{x}-X_{y}$ gives sense of direction for the system. $W>0$ iff $Y_{x}/X_{y}>1$.
Since $\frac{Y_{x}}{X_{y}}=\rho \frac{d Y}{d X}$ and $= Y \frac{d}{d X} \rho + \rho^{2} > 1$. A sufficient condition for this is $Y \frac{d}{d X} \rho > 0$, assuming $\rho^{2} > 1$. If $\rho$ is odd, then $\frac{d}{d X} \rho$ is even and if $X < 0$ implies $\rho X=Y > 0$ then $\rho < 0$, if $X > 0$ implies $\rho X =Y > 0$ then $\rho > 0$. If $\rho$ is even, then $\frac{d}{d X} \rho$ is odd and $X < 0$ implies $\rho > 0$, $X> 0$ implies $\rho > 0$. When $\frac{d}{d X} \rho=\frac{d}{d X_{1}} \rho$, we assume $Y,Y_{1}$ are of the same sign. A mirror, corresponding to an irreducible component in a first surface, is given by $\rho= c$, for a constant $c$ and the second mirror is given by $(X_{2},Y_{2})$ where $\rho=c$. For the mapping $\gamma \rightarrow \zeta$, we must consider the restriction to lines.

Assume $\gamma \rightarrow P_{0}$, a point on the boundary, is given by $d F (\gamma)=0$ iff $< d F, \gamma >=0$. In the spiral case, we are assuming $dF (\gamma) \neq 0$ and we can assume the spiral has representation with measure that is not reduced. Consider now the intersection between the curve $\gamma$ and the normal $N$.  We assume the intersection locally algebraic, but not necessarily locally injective. For $N=\{ P(d \gamma)=0 \}$ and $\gamma$ such that $dF(\gamma) + P(d \gamma)=0$, assume $RP(d \gamma) \rightarrow s d \gamma = \tilde{d \gamma}$ continuous, where $R$ is assumed locally algebraic. If $R$ is absolute continuous, we have that when $\tilde{d \gamma} \neq 0$ and $P( \tilde{d \gamma}) \neq 0$, then $P(d \gamma) \neq 0$. In the transversal case, we assume the measure is reduced on irreducible components.

The second problem we should consider is to determine the geometric sets, that are given by a fixed $\gamma_{1}$ (and $\gamma_{2})$ as above, given that $(I)$ can be continued to these systems.
In the case of two mirrors, we have not excluded spirals in the domain for $(\tilde{I})$.

\subsubsection*{The transmission property}

In connection with what we denote the transmission property (TP) in this article, $u(\beta) \sim u(\beta^{*-1})$, we are discussing the class of symbols for which $u+c(\beta) \sim (u+c)(\beta^{*-1})$, without assuming the same property for $u$.
A sufficient condition for $\{ f-c \} \sim \{ f-0 \}$ is isolated singularities for instance $f \in L^{1}$. The same condition can be given for $\phi= \log f \in L^{1}$ and $\phi + c (\beta)=\phi(\beta_{1})$, for some $\beta_{1}$ in the domain for $\phi \in L^{1}$, corresponding to an interpolation property. The advantage with using several ($>1$) mirrors in a localization of the symbol is possible surjectivity, which implies injectivity for the transposed action (\cite{Martineau}). 
Assuming surjectivity, of $f+c \in (I)$, we can find $F$ such that $F(\tilde{\gamma})=f+c$. Let $T_{1}F(\gamma)=F(\tilde{\gamma})$, corresponding to continuation to the lower half plane in the phase, where $T_{1}$ is assumed continuous at the boundary. Further, at the boundary $T_{1}F(\gamma)=\zeta F(\gamma)=F(\tilde{\gamma})$ and $\mid \zeta \mid=1$. In a Puiseux expansion,
$\chi (=y/x)) =a_{0}+a_{1} t +\ldots$, where $t = \sqrt{\zeta}$ and $a_{1} \neq 0$ means that the order for the critical surface is $1$. Note that if $\beta . \beta^{*-1}=1$ then $F(\beta)=\zeta F(\beta^{*-1})$ and $\mid \zeta \mid=1$, assuming that $F$ has TP. In a desingularization, it is to reach $V$ in the domain for $F(\tilde{\gamma})$, sufficient that $F$ is analytic over the final segment and the initial one, over the other segments $F$ can be chosen in $H'$. That is if the path is divided into $N$ segments, $V_{N}$ must be in the domain for analyticity for $(X_{N},Y_{N})$ and $F(\tilde{\gamma}_{N})$. The mapping $\gamma_{T} \rightarrow \zeta_{T}$ in $F(\ram \gamma)(\zeta)=F(\gamma)(\zeta_{T})$ is locally 1-1 in the reflection model, why we must have a multivalued image. Note that $V$ is not in the domain for $F(\gamma)$. Consider now $\frac{u(\tilde{\gamma})}{u(\gamma)}=\frac{{}^{t} \tau u(\gamma)}{u(\gamma)}$. Then, $u$ is odd
iff ${}^{t} \tau u(\gamma)/u(\gamma)$ is odd. We are assuming $u$ odd on a small disk and discuss if the property is preserved under onesided continuations. If $u$ has TP onesided, then $T_{1}$ is not equivalent with $T_{2}$. That is $T_{1}f(z)=f(\alpha.z)$, where $\log F(\alpha.z) \sim \log F(\alpha^{*-1}.z)$ and $T_{2}f(z)=F(z.\gamma)$, $\log F(z.\gamma) \sim \log F(z.\gamma^{*-1})$, but we do not assume $\big[ I,T_{1} \big]=\big[ T_{2},I \big]$. Note also the difference between continuation over a removable set and continuous continuation over larger distances. If $F$ has TP that is preserved under onesided continuations,
$\big[ T_{1},\tilde{I} \big] F \in H$ iff $T_{1}F \in H$. We construct the inverse according to Julia (\cite{Julia24}) using for instance Abel's equation and $\Phi(\alpha.z) \rightarrow \zeta + \eta_{+}$ and $\Phi(z.\gamma) \rightarrow \zeta + \eta_{-}$  The domain $D$ for $F$ and TP can, given the systems be determined in this way. If the domain can be continued using $T_{1},T_{2}$, then the mirrors are seen as symmetry axes for the domain and the property TP. Thus, if $z \in D$ then $z^{*-1} \in D$ and conversely. If $\alpha.z \in D$ then $z^{*-1}.\alpha^{*-1} \in D$ and conversely. Thus $\big[ \tilde{I},T_{1} \big] \sim \big[ T_{2},\tilde{I} \big]$ and $D$ is still symmetric for the continued curve. If the distance between two mirrors is small, then this corresponds to a symmetric disk. A bigger distance corresponds to a TP property for the solution, that is that $V$ corresponds to solutions with this property. If $D$ is symmetric with respect to two axes, we have a twosided continuation. If we consider the mirrors $R_{0},R_{1}$ as irreducibles in a first surface, the boundary can be seen as analytic. The domain has TP if $\tilde{I}F=F\tilde{I}$, but if we have symmetry with respect to $T_{1}$, we do not necessarily have symmetry with respect to $T_{2}$. Note that the problem to reach $V$ using the model, does not require TP, but the TP property means that we can reach $V$ through transversals.

\vsp

If $\Omega$ is the domain for $F(\gamma)$, we are assuming that $V$ is in the complement to $\Omega$. Assume $F(\tilde{\gamma})$ is the continuation through $R_{0},R_{1}$, then we have $V \subset \tilde{D}$, the domain for $F(\tilde{\gamma})$. Thus the symmetry property TP is relative the complement to the domain $\Omega$. Note that we are discussing the domain $\tilde{D}$ locally close to the segment $\gamma_{N}$. Note also that the symmetry property can be seen in three steps,
$\tilde{I}F(\gamma)(\zeta) \sim F(\tilde{I} \gamma)(\zeta) \sim F(\gamma)(\tilde{I} \zeta)$.
The original problem that was suggested to me by Prof. Vuorinen at Helsinki University,
was to describe the clustersets when the complement set has a line with three points. We assume in this article two reflection points, taken at (possibly different) points in time. The third point is the origin
which is necessary to determine how the mirror is oriented.

Assume $N_{j}=Y_{j}/X_{j}$ and $F(X,Y) \sim F_{1}(N)$. Thus, $N_{j}=const$ on $R_{j}(X_{j},Y_{j})$.and $\{ F_{1}=const \}$ corresponds to $\{ f=c \}$. We are assuming $F_{1}$ maps constants onto constants, where $F_{1} \in H'$. If $g_{1},g_{2}$ are algebraic (continuous, AC) then $g_{1}g_{2}$ is algebraic (continuous, AC). If $D$ is the domain for analyticity for $F(\gamma_{N})$, consider $D'$ the domain for absolute continuity (AC).
If $V \subset D'$, we are done, If $D''$ is the domain where $F$ maps $\{ N=c \}$ onto $\{ f=c \}$ and if $V \subset D''$, we are done. Note that since in the problem the systems $(X_{N},Y_{N})$ can be varied, we can select $(P,Q)$ such that $D'' \supset V$. 

\vsp

If we assume $F \in \mathcal{D}_{L^{1}}'$ over $(X_{j},Y_{j})$ analytic (normal) with $\log F(\gamma_{j}) \in L^{1}$, we have singularities of finite order, for all $j$. In the case when $F(\gamma_{j})$ analytic, we can assume $V_{N} \rightarrow V_{N+1}$ is continuous and $V=\lim_{j} V_{j}$, where the limes is over surfaces of finite order. Given TP reflection using $R_{0},R_{1}$, this is transposed to a normal model and we have singularities of finite order. If we assume the generator for $L^{1}$ taken in phase space $\phi_{0}$, the $K_{0}=(\mbox{ ker }\phi_{0})$ leaves space for spiral approximations. If we assume $K_{0}$ trivial (regularizing) and that we do not have a normal model, we can not have TP!

Note also if we assume in $\gamma_{1} \rightarrow \ldots \gamma_{N}$, where $F$ is analytic over the first and last segments, we have that $V_{1},V_{N}$ of finite order, but the geometric sets corresponding to the trajectories between the systems $(X_{j},Y_{j})$, given by a contact transform, are possibly of infinite order or undetermined. $V$ is considered as the limit also in this case. In the case when $T$ denotes time, we notice the difference between $F(\tilde{\gamma}_{T})$ and $F(\tilde{\gamma})_{T}$. 

Assume $\log F(\gamma)=\Sigma \phi_{0}(\zeta + \eta_{j})$ over the complement to $K_{0}$. Over the same set $F(\gamma) \sim \Pi_{1}^{\infty} e^{\phi_{0}}(\zeta + \eta_{j})$. If we assume $\phi_{0}$ such that $e^{\phi_{0}}$ analytic, we have $\phi_{0}(\zeta + \eta_{j}) \rightarrow \Omega_{j}$ continuous and $\Omega=\cup_{1}^{\infty} \Omega_{j}$ (not necessarily disjoint). The mirrors are considered related to different systems, why if $\zeta + \eta_{j}$ generates $R_{0}$ and $\zeta + \tau_{j}$ generates $R_{1}$, the different leaves are defined by $\gamma$.
Note also that $F$ can be reduced as a measure, without being analytic.

\vsp

Note the following example. Assume $(I)$ paths in the dynamical system and that $1/\beta \in (I)$, then we have existence of $\tilde{\beta} \in (I)$, such that $u(\frac{1}{\beta}) = u(\tilde{\beta})  =\Sigma_{j} u(\beta)(\zeta + \eta_{j})=\Sigma_{j} u((\ram)_{j} \beta)(\zeta)$. Note also that if $u$ is not algebraic in $\beta$, then $u$ can still be algebraic in $1/\beta$ and in $\tilde{\beta} \backslash \beta$. Assume $v(\beta)=\phi_{0}(\frac{1}{\beta}) - \phi_{0}(\beta)$, then $v(+ \infty)=-v(+0)$ and $v(-\infty)=-v(0)$. In this case $v(\frac{1}{\beta}) \sim -v(\beta)$.

\subsubsection*{Algebraic inverse}
Assume $\varphi$ defines $f(z) \rightarrow z$. If $id=I$ acts on $(I) \ni f$ is defined through some $g \in \mathcal{G}$ (reflection through the origin) we must discuss $\big[\varphi,I\big]-\big[I,\varphi\big]$. We are assuming $\big[\varphi,I\big](f)=-z$ and $\big[I,\varphi\big](f)=\varphi(-f)(z)$. If $f(-z)=-f(z)$, then $\big[\varphi,I\big]=\big[I,\varphi\big]$, that is $f$ is odd. More generally, consider the restriction $\big[R,\varphi\big](f)=z \mid_{\Omega}$ defined as $z$ for $z \in \Omega$ and $0$ for $z \notin \Omega$, where $\Omega=\mbox{ supp }f$. Thus, $\big[\varphi,R\big](f)=(f \mid_{\Omega}) \rightarrow z$, that is on $\Omega=\mbox{ supp }f$, $R$ can be defined algebraic.
Assume for a homomorphism $g$, $\big[I,\varphi\big](f)=\varphi(g f)$ and $\big[\varphi,I\big](f)=\varphi(f)(gz)$,
that is if $g f = f g$, we have $\big[I,\varphi\big]=\big[\varphi,I\big]$. Assume that $f$ is not odd, then using reflection, we have existence of $u$ such that $f(u)=-f(x)$. Extend the domain associated to $\tilde{f}$ to include $U$. Then if $\tilde{f}$ is analytic, we have $z \rightarrow u$ is continuous. Assume $w$ a bijection,
such that $\varphi(g f)=\varphi(f w g)$. Thus, $\big[I,\varphi\big] \simeq \big[\varphi,I\big]$.

\subsubsection*{Irreducibles}
Let $L=L_{1} \cup L_{2}$ implies $L_{1}=\{ 0 \}$. Alternatively, assume $L_{1}$ is removable (algebraic)
that is $f \in (I)(L)$ and $f=f_{1}f_{2}$ implies $f_{1}$ algebraic. Assume $\mu_{L}$ a reduced measure
associated to an irreducible $L$ and let $\tilde{\mu}_{L}$ be the result of an algebraic continuation.
If $\tilde{L}$ is  the topological continuation associated to $\tilde{\mu}_{L}$, we have that this can be selected 
irreducible (modulo algebraic sets). Note that we do not assume that irreducibles contain the $\infty$. 
Alternatively, we could let the mirrors depend on a parameter $\mu$ (time) and let the reflection
be at different times. In this case, we need an origin as a third point to orient the plane in space. If the reflection curve is on one side of one mirror, we do not have spirals, but if it is
on one side of two mirrors, we must discuss spiral approximations. If we consider $\psi$
as $\psi_{1}\psi_{2}$, where $y=\psi_{2}x$ is onesided with respect to the mirror $R_{1}$ and $\psi_{1}y$
is onesided with respect to the mirror $R_{2}$, then if $\psi$ is irreducible, we have that $R_{1}=R_{2}$. When the TP is present, so that we can transform the reflecting curve to a transversal,
the result of two reflecting mirrors is an irreducible.

\subsubsection*{Reflection and onesidedness}
Assume $L=\alpha \beta \gamma$ describes a reflection using two mirrors $R_{1},R_{2}$ and after
moving the reflection to the origin and $R_{1}$ to the imaginary axes using rotation and translation, let 
$\beta \rightarrow \beta^{*}$ be reflection through the origin and $\beta^{* -1}$ the result of
reflection through the origin with reversed direction of path. Let $\tilde{L}=\alpha \beta^{* -1} \gamma$.
If $f$ has the property that $f(L)$ analytic iff $f(\tilde{L})$ analytic, the the reflection can be compared
with a transversal approximation, In this case $\tilde{L}$ can be regarded as on one side of a plane
and $R_{1},R_{2}$ are first surfaces with reflection points as singularities. Let us write the reflection
using two mappings $T_{1}f(z)=f(\alpha z)$ and $T_{2}f(z)=f(z \gamma)$. If for example
$T_{1}f(\beta^{* -1})=-T_{1}f(\beta)$, then obviously the analyticity is preserved.  The properties of the inverse mapping $T_{1}f(\beta) \rightarrow \beta$
is dependent on the monodromy condition.

The condition that $f$ is odd
is written as $\big[ \delta,f \big] = \big[f,\delta \big]$. We assume that for any function $f$,
locally analytic in $z_{0}$, there is a disk with center $z_{0}$ where $f$ is odd. We can assume
that $f$ is analytic independent of orientation for boundary to the disk. We will continue the disk using the mappings
$T_{1},T_{2}$. Note that $\big[ I , T_{1} \big](f)=I(f)(\alpha \beta)=f(\beta^{* -1} \alpha^{* -1})$
and if we let $\gamma \sim \alpha^{* -1}$, $=\big[ T_{2},I \big](f)$. Let $\mathcal{G}_{1}$ be 
continuations one sided to the lower halfplane. If we assume $f$ algebraic in $\beta$ in the disk,
we can assume the complement in a domain of holomorphy, that is $\alpha$ has an analytic parameterization. The argument can be repeated to $\mathcal{G}_{2}$ and continuations to the upper halfplane. Note that an (locally) algebraic extension assumes an oriented boundary of the disk.  

Finally note that global analysis means that our propositions are invariant for change in
local coordinates. When time is considered as a variable, this means a transportation in time.
Thus the reflection model as we have described it, is not global in the sense above, but it can with
advantage be combined with global analysis in localizing the symbol.

\subsubsection*{Schur's lemma}
Schur's lemma gives that if $A$ is a linear mapping $N \rightarrow M$, where $N,M$ are finite dimensional and $w_{1},w_{2}$ irreducibles
such that $A w_{1} = w_{2} A$, then $w_{1} \sim w_{2}$ or $A=0$.  Another version of the lemma is with the same conditions, we have existence of $A^{-1}$ or $A=0$. The Schur's lemma can be generalized to infinitely dimensional spaces. We further note that the number of linearly independent representations is the same as the number of irreducible representations. 

\vsp

Assume $J$ an irreducible ideal, that is if $J=J_{1}J_{2}$ and $N(J)=N(J_{1}) \cup N(J_{2})$, where $J_{i}$ analytic, irreducible and such that $N(J_{i})=\{ 0 \}$, for one $i$. Thus, if $f \in J$ and $f=gh$ implies $g=id$ ($\sim \delta_{0}$). Thus, $f=\big[ I,h \big]$ or $\big[ g, I \big]$ and $\big[I,h \big]=\big[ g,I \big]$ implies $I=0$ or we have existence $I^{-1}$ and $h \sim g$. Assume now $g \in \mathcal{G}$
gives analytic continuation, that is $\tilde{f}=g f$, that is $\mathcal{G}$ is acting outside the range.
We consider $(I) \rightarrow (X) \rightarrow \Omega$ and parallell $(J) \rightarrow (W) \rightarrow V$.
For $(J)$ to be an analytic ideal, it is sufficient to have the lifting principle, for instance $W$ algebraic.
 In the scheme above, if $F(\tilde{X})=G(W)$, over the set where we have existence of $G^{-1}$, we have $W=G^{-1}F(\tilde{X})$. Particularly, if $\psi(\gamma)=\overline{\gamma}$, that is complex conjugation, the proposition is existence of $W$ a dynamical system with $\sigma=\{ \gamma \quad W_{1}/W_{2}=const \}$, such that $\varphi(\psi(\gamma)) \in V$, the geometric event we want to predict. If we assume the continuation according to Lie, that is $X_{1}/X_{2}=W_{1}/W_{2}$, then if  $\gamma \in \{ X_{1}/X_{2}=const \}$, then we have $\psi(\gamma) \in \sigma$.

Given $f$ an entire function with $\{ f-c \}$ irreducible, we refer to this set as a first surface.
Assume $\{ S_{j} \}$ first surfaces, and $S_{j} \rightarrow S$,$S_{j} \rightarrow S'$ implies $S \sim S'$ that is the limits are conjugated. For a regular first surface, we have that $S=S'$. For an entire function almost every first surface is regular (\cite{Nishino68}). Thus, assuming a normal model
where all singularities can be considered in first surfaces, we must have $\big[ \varphi,I \big] \sim \big[ I, \varphi \big]$ locally, close to the constant surfaces and where reflection is adapted to the point of intersection with the transversal. Thus, the limit is not dependent on from which side we do the approximation.

\newtheorem{prop4}[prop1]{Lemma}
\begin{prop4}
 For an entire symbol, where all singularities can be approximated using a normal model, then the inverse mapping $\varphi$ is algebraic in the sense that it commutes with $I$ locally, close to the first surfaces.
\end{prop4}

Note that if $\log \mid f \mid \in L^{1}$ and $g \in \mathcal{G}$ implies that $\log \mid g f \mid \in L^{1}$. This means that the domain for the linear $gI (=A)$ is standard complexified. Further $f$ will have finite order singularities and the order for a singular point is defined by the order for the point on any complex line.  For instance, if the reflection is with reflection point the origin, then real axes is normal to the imaginary axis as mirror and the reflection result is complex conjugation. If $V$ is parabolic, we have that we can use monotropy in the reflection, and assume the imaginary axes holds all reflection points. If we use two mirrors and one reflection point (origin), we have $g_{2}g_{1}X=\psi(X)$ and $g_{3}g_{2}^{-1}\psi(X)=g_{3}g_{1}(X)$, in this example the two mirrors are the axes.

Define ${}^{t} \psi(X) \tilde{F}(X)=\tilde{F}(\psi X)$ by $\tilde{F}(X) \tilde{F}(\psi X) \sim \mid \tilde{F}(X) \mid^{2}$ (\cite{Riesz_23}) such that $\log \tilde{F}(\psi X) = - \log \tilde{F}(X)$. Invertability assume that $u \rightarrow -u$ is surjective and that the path $X \rightarrow \psi X$ is an analytic continuation (Thullen). Note that if $\psi$ corresponds to a conformal reflection, then analyticity is preserved but not microlocal properties.  Note the example $f=e^{\phi}$ with $\phi \in L^{1}$ $\Leftrightarrow - \phi \in L^{1}$. Thus, $f$ and $1/f$ have finite order singularities, but this does not imply that $1/f$ is analytic. However if $X \rightarrow \psi X$ is a path for analytic continuation, we have that $1/f$ must be analytic.

In the scheme $(I) \rightarrow X \rightarrow \Omega$ and parallell $(J) \rightarrow W \rightarrow V$, the mapping $(I) \rightarrow (J)$ corresponds to ``blow-up'' if we do not change origin and corresponds to a simulation in this case, that is we use $(I) \rightarrow X \rightarrow W \rightarrow (J)$. If we change origin (of time), that is we use $(I) \rightarrow \Omega \rightarrow V \rightarrow (J)$ then $(J)$ is a true geometric ideal. The mapping $(I) \rightarrow (W) \rightarrow V$ does not involve ``blow-up''.

Assume $X \rightarrow \psi X$ is a reflection such that $\psi X$ is harmonic conjugation.
Then $\psi$ is pure if $\psi(X + i \psi X)=-i(X + i \psi X)$. 
If $u , \psi (u)$ are closed forms, we have that $u$ is pure, that is $u = a d z$ with $a$ locally analytic.
If $u$ is closed, $\psi (u)=-i u$, then $u$ is analytic and $X \rightarrow \psi (X)$ is permitted for analytic continuation (\cite{AhlforsSario60}). 
A transversal is a line that cuts two other lines. In the reflection model, this is the case when we have two mirrors and TP.  Consider the following example. Starting with $\{ \phi - \mu \}$, where $\mu$ is constant, the transversals are given by $q \bot \mu I$, where $q$ are quasi orthogonal. Thus, in the case with one mirror, we assume the reflection points are given by zero's to $q$ quasi orthogonal with real coefficients. If we are considering the case with two mirrors, we are assuming the case with complex zero's on the axes. Particularly, if the resolution is given by $q_{j}$ such that $\psi q_{j}=q_{j} \psi$, the support is symmetric with respect to $\psi$ (cf. standard complexification).

\subsubsection*{Geometric reflection and harmonic conjugation}
We define for $u=a d x + b  d y$,$u^{*}=-b d x + a d y$.
Define $\tilde{u}(z)=u(\overline{z})$, then we have if $u$ is even, this implies $u^{*}(z)=\psi(u)$
is odd. When $u^{*}$ is odd, $u^{*}(\psi(z))=\psi(u^{*}(z))$. For an algebraic $H$, we have $H(e^{-\phi})=1/H(e^{\phi})$, that is $\widehat{H}(-\phi)=1/\widehat{H}(\phi)$, or $\log H(e^{-\phi})=-\log H(e^{\phi})$, that is $H$ is odd in the phase.  The condition $u(\psi(z))=u(z)$ does not exclude a zeroset that contains spiral sets. However, for $u$ even, we have existence of $u^{*}(z) \rightarrow z$
continuous. Thus $\psi$ does not preserve algebraicity or invertibility.
Further, if $u^{*}$ is analytic and $u$ closed, then $u$ is analytic, why we have existence of $\gamma$ a path for analytic continuation.

Consider the following example.
$\overline{u}(x,-y)=a d x - b d y$ and $\overline{u}^{*}(x,-y)=b d x + a d y$. Further, if we let
$w(z)=a(z) d x + b(z) d y$ and $w^{\sim}(z)=a(\overline{z}) d x - b(\overline{z}) d y$, then
$w^{* \sim}(z)=-b(\overline{z}) d x - a(\overline{z}) d y = -w^{\sim *}$, that is if $w$ is even, then
$w^{*}$ is odd. Assume $\psi \in \mathcal{G}$ is reflection through the origin, then $u \psi= \psi u$ implies $-u(x)=u(-x)$ and $u \psi=u$ implies
$u(-x)=u(x)$. Assume $\mathcal{G}_{1}$ corresponds to geometric reflection
over the real axes, $u^{\sim}$ and $\mathcal{G}_{2}$ corresponds to harmonic conjugation $u^{*}$,
with the property $u^{* \sim}=-u^{\tilde *}$. Then we have for example $g_{1}g_{2}=-g_{2}g_{1}$.
An arbitrary form is not necessarily odd or even, however if $u^{*}$ is odd, then $u(-\overline{z})
=u(\overline{z})$ and if $a,b$ are even, then $u$ is even (\cite{AhlforsSario60}). Assume the representation of the phase $\phi$
is such that we can find a parameter $\mu$ so that $\phi - \mu$ is odd. For these phases $d \phi$ is even
and $(\phi - \mu)^{*}$ is even. Further, $\phi_{\mu}^{*}/\phi_{\mu}$ is odd in argument.

\vsp

If $\mathcal{G}$ defines analytic continuation, we note the concept of monogenity which means that 
every element $gf$, for $g \in \mathcal{G}$, can be defined from any other element $hf$ with $h \in \mathcal{G}$
for instance $\tilde{f}=gf$ and $\tilde{\tilde{f}}=v f$ implies existence of $w \in \mathcal{G}$
such that $\tilde{\tilde{f}}=wg f$ and this corresponds to an interpolation property. Note that if
$v=wg$, then  this defines a set of invariance in $\mathcal{G}$, that is it is necessary to discuss
disjoint supports. In particular, we note that if $g \in \mathcal{G}$
has a representation with a reduced measure in $\mathcal{E}^{' (0)}$, we do not have the center case and 
we have existence of an inverse $g^{-1} \in \mathcal{G}$. Thus the concept of hypoellipticity is useful
in connection of analytic continuation. When $\mathcal{G}$ is geometric reflection, the condition
that the path is on one side of a mirror excludes the center case (cf ArXiv).

\subsubsection*{Involution and the transmission property}
Assume $\gamma_{0}$ corresponds to $(X_{0},Y_{0})$ and involution $\Sigma_{0}$ and $\gamma_{1}$ corresponds to $(X_{1},Y_{1})$ and $\Sigma_{1}$ and that $\Sigma_{0} \rightarrow \Sigma_{1}$ continuous. If $Y_{0}=\frac{-d \phi_{0}}{d x}$ and $X_{0}=\frac{d \phi_{0}}{d y}$, we can write
$\frac{Y_{0}}{X_{0}}=\frac{-d \phi_{0} / d x}{d \phi_{0} / d y}$ on 
$\Sigma_{0}$. The discussion is now about contraction of formulas. 
We define $\Delta$ and $\bot$ as independent of $\Sigma_{j}$. Assume $\phi_{0}^{\delta}$ defined in a neighborhood of $\Sigma_{0}$, such that we have existence of inverse $\phi_{0}^{\delta} \rightarrow \gamma_{T}$ over $\Sigma_{0}$.
The symmetry condition is $\frac{d y / d t}{d x / d t}=\frac{- d \phi^{\delta} / d x}{d \phi^{\delta} / d y}$ 
on $\Sigma$. The regularity conditions imply that 
$\frac{d \gamma}{d t}=\frac{d \gamma}{ d x}\frac{d x}{d t} + \frac{d \gamma}{d y}\frac{d y}{d t}=0$ 
implies $t=0$. The symmetry conditions means that the regularity conditions 
implies that the lifting limit  is $=id$. 
When the vorticity $W=0$, the symmetry condition means that $\frac{d^{2} \gamma}{d x^{2}}+\frac{d^{2} \gamma}{d y^{2}}=0$. 

 If $\phi \in L^{1}(D)$, for a domain $D$, we have isolated singularities $S(\phi)$ and for $f=e^{\phi}$, we have finite order singularities (\cite{Riesz15}). 
Assume $V$ is given by a polygonal between $N(\phi_{0})$ and $V$. If $\gamma=e^{\phi}$ and 
$\phi \neq - \infty$, then $\frac{d \gamma}{d T}=0$ iff $\frac{d \phi}{d T} =0$, that is the system can 
be referred to the phases and the mapping $D/S(f) \rightarrow D/S(\phi)$ can be assumed algebraic. The mirrors correspond 
to spectrum relative the dynamical system, that is $Y/X=const$.
If $(X,Y)$ are analytic, then $I(V)$ is defined using algebraic topology, but this property is not preserved 
when we use $(X,Y) \rightarrow (X_{1},Y_{1})$ continuous, where the range is analytic. 
If in the mirror model we have a path with self intersections, we get 
inconsistencies in the orientation of path. In this case we can consider a covering of leafs.  An example when self intersections do not cause inconsistencies of any importance, is
the mirror points correspond to different points in time.

If we assume $f$ exponential, then the continuation is completely determined by the phase.
Assume for instance $f=e^{\phi}$ and $g=\beta e^{\phi}$, where $\phi \in L^{1}$. On the support for $\phi$, 
we have $g=0$ iff $\beta=0$ and $\frac{d g}{d x}=\frac{d \beta}{d x} e^{\phi}=0$ iff $\frac{d \beta }{d x}=0$
and $\frac{d}{d x}g=0$ iff $\frac{d}{d x} \log \beta=-\frac{d \phi}{d x}$ iff $\frac{d g}{d x}=0$.

We assume for the paths, the property that a boundary point can be reached through a polygon using mirrors and blowup. 
We assume for the translations $<\eta, \tau_{h}x>=<\tau_{h}' \eta,x>$ with $\tau_{h}'$ locally algebraic.  Assume $\phi \in (\mbox{ ker }h)$ and $h^{2}=1$,$h^{*}=h$. Assume $\phi(x)=<\eta,x>=N(x)$ and $h(\phi)(x)=[H,N](x)$ and $[H,H]=1$. If $C=HN-NH$ then $HHN - NHH=HC + CH=0$.
The condition $[H,H]=0$ on $V$ can be read as $H \bot N$ and if $H^{2} \bot N$ on $V$, we have $N=0$ on $V$. Assume the mapping between the dynamical systems is linear and continuous and the quotients are preserved by the continuations, within a constant and that the orientations are preserved.
Using Wiener's argument, since $f \in L^{1}$ iff $\mid f \mid \in L^{1}$, we can assume movement is translation and rotation. Assume $L^{1} \ni \phi \rightarrow \phi_{0} \rightarrow V$, where $\phi_{0} \in L^{1}$ is generator, that is $\phi$ has representation on $\mbox{ supp }\phi_{0}$. For $\phi_{0}$ we have regular approximations that can be chosen as normal. For a normal model we have $\phi_{0} \in (\mbox{ ker }h)$ with $h^{2}=1$ and $h^{*}=h$. For $(I)=(\mbox{ ker }h)(\Omega)$, where $\Omega$ is pseudoconvex, we have that tangents $\sim$ normals.

\vsp

In the reflection model $\gamma_{0} \rightarrow \gamma_{1} \rightarrow \ldots \rightarrow \gamma_{N}$ 
with respective geometric sets $V_{0},V_{1},\ldots,V_{N}$, if the problem is to localize $V$,  we do not require $V_{0} \rightarrow V_{1} \ldots \rightarrow V_{N}$ 
continuous. Assume $F^{-1} \sim P(D) \delta_{0}$ (in $\mathcal{D}_{L^{1}}'$)
locally at the boundary. If the representation is extended to the entire 
space, we can give the solution exactly in $\mathcal{D}_{L^{1}}'$ (and in 
$H'$). With this approach, denote $F_{\delta}$ for an approximation of $F$ close to the boundary and we have that $F_{\delta}$ can be represented as a parametrix to 
$F^{-1}$. If $\Gamma$ is very regular we have that $P(D)$ is hypoelliptic 
implies that $F_{\delta}$ is very regular. 

Consider $F(X,Y) \rightarrow (X,Y)$. If we have existence of $F_{\delta}$, then the problem $F_{\delta} \rightarrow \delta_{0}$ regular corresponds to local existence of $F^{-1}$. Assume $G$ satisfies $\frac{d F}{d x} \frac{d G}{d y}=\frac{d F}{d y}\frac{d G}{d x}$ where $\frac{d F}{d x}=Y$ and $X/Y=-\frac{d G / d x}{d G / d y}$, that is we consider continuations $G$, that preserves the quotient of right hand sides, that is if the involutive set is $\Sigma$, then we consider $I(\Sigma)=\{ G \quad Y/X=\frac{d G / d x}{d G / d y} \}$

\subsubsection*{Discussion}

Concerning Poincar\' e's objection (\cite{Poincare87},253bis) about continuations according to $\tilde{Y}/\tilde{X}=Y/X$,
with $\tilde{X}=Z X$ and $\tilde{Y}=Z Y$ and $\frac{d t_{1}}{d t}=Z$. If $M_{0} \rightarrow M$ ($dt$) at $t=\tau$ and $M_{0} \rightarrow M'$ ($dt_{1}$) at $t_{1}=\tau$. Further, $F_{0} \rightarrow F$ at 
$t=\tau$ and $M_{0} \rightarrow F_{0}$ at $t=0$ and $F_{0} \rightarrow F'$ by $t=\tau$, then there 
are not necessarily $\sigma$, such that $F=F'$ at $t=\sigma$. In the reflection model, we consider time as a local parameter. The contact transforms between the systems are not assumed to preserve invariants.
The characteristic sets as lineality and orthogonality are formed in the blow-up to the local ideals
corresponding to the systems. We prefer to use the parameter $T$ as a global parameter and $t$ for local parameters. Note that if $F(\tilde{\gamma_{T}})=\tilde{F}_{T}(\gamma)$ and $G_{T}$ is the lifting function to the event $V$, then we are not assuming $\tilde{F}_{T}=G_{T}$. We are not assuming $\tilde{F}_{T}$ analytic, only that the mapping $\tilde{\gamma}_{T} \rightarrow \zeta_{T}$, locally defined, is continuous.

\subsection*{ Representation at the boundary }

Note that in a discussion of sections, we require completeness, that is dependence on all variables
in space. Also note that
a sufficient condition for $\log X_{T} \rightarrow V$ to be continuous, is that the logarithm is
analytic.  We can form a domain, as the largest for which the mapping is continuous.
Finally, note that the ideal over singularities given by $\frac{d}{d T} \log f_{} = const.$ must be radical. 

\vsp

Assume $N(e^{\psi})=e^{\tilde{N}(\psi)}=e^{<\eta_{T},\psi>}$ and if $N$ is algebraic, we have that $< \eta_{T},e^{\psi}>=1$ iff $e^{<\eta_{T},\psi>}=1$ and $\widehat{I}N(\psi)=N\widehat{I}(\psi)=\widehat{N}I(\psi)$. If $\mid < \eta_{T},\psi> \mid \leq C_{T} \mid \psi \mid$, with $C_{T} \rightarrow 0$, as $T \rightarrow \infty$, then we have that $N$ can be represented by a measure of finite type.
If $N$ is algebraic, we have $\big[ \widehat{I}, N \big](\psi)=\big[ I, \widehat{N} \big](\psi)$ and $I \sim \delta_{0}$. If we let $H_{1}(\Omega)$ be holomorphic and non-constant functions on $\Omega$, then  if $e^{\psi} \in H_{1}(\Omega)$, then $e^{\psi+c} \in H(\Omega)$.
If $\mid C_{T} \mid \rightarrow 0$, as $T \rightarrow \infty$, we have $\mid \eta_{T}-P(1/T) \mid < \epsilon$, as $T \rightarrow \infty$. In this manner we can compare $\eta_{T}$ with Legendre in the infinity, modulo monotropy. Note also that $\mid N(e^{\psi}) \mid \leq e^{\epsilon \mid \psi \mid}$, as $T \rightarrow \infty$, that is $N$ has type zero with respect to $\psi$. We are assuming $N$ is acting only in phase space, that is $e^{<\eta_{1}+\eta_{2},\psi>}f=N_{1}(\psi)N_{2}(\psi)f$ and $e^{\Sigma_{1}^{\infty}<n^{j}_{T},\psi>}f=\Pi_{1}^{\infty}N_{j}(\psi)f$. For convergence, it is necessary that $\mid \eta_{T}^{j} \mid \rightarrow 0$, as  $j \uparrow \infty$ and $T$ fix.

\vsp

Assume $<\eta,\phi>=\eta(\phi)$ for $\phi \in H(\Omega)$, we have existence of a constant coefficients polynomial differential operator $P(D)$, such that $P(D) \eta \sim \delta_{0}$
(modulo regularizating action), then $\eta$ has representation with one sided support. Assume we have
existence of an inverse $\eta^{-1}$ with respect to convolution, then we are assuming $\widehat{\eta^{-1}}$
is a polynomial. Note that if $f,g$ are hypoelliptic symbols and $\int_{\Gamma} fg d \sigma=0$
then either $\Gamma$ is algebraic (removable) or $fg \equiv 0$ on $\Gamma$, which must be a bounded set.
 
\vsp

If $\varphi \geq 0$ imples $\overline{\varphi} \leq 0$ and that $\Omega=\{ \mbox{ Im }e^{\varphi} < \lambda \}$ unbounded, that is $\Omega \sim \{ e^{\varphi} < \lambda+e^{\overline{\phi}} \}$. Then, when $e^{\varphi}$ is parabolic, we have that $\{ e^{\varphi} < \lambda' \}$ is unbounded.

\subsubsection*{Hypoelliptic differential operators}

Consider now the condition $\mbox{ Im }P(\xi) / \mbox{ Re }P(\xi) \rightarrow 0$, as $\mid \xi \mid \rightarrow \infty$, which we denote $\mbox{ Im }P \bot \mbox{ Re }P$.
Assume $\mbox{ Re }f$ hypoelliptic and $f$ normal, then we have that $\mbox{ Im }f \bot \mbox{ Re }f$. Further we have that
if $\mbox{ Im }f$ bounded, we have $\mbox{ Im }f \bot \mbox{ Re }f$. Further $\mbox{ Im }f$ unbounded does not mean
that $\mbox{ Im }f \bot \mbox{ Re }f$. 

\vsp

If $P$ is hypoelliptic implies $NP$ is hypoelliptic, then $N$ has algebraic definition. If we have existence of a segment in the infinity, such that $NP/P \equiv 1$, then we have that $P$ is not hypoelliptic, in the infinity. The condition $c_{1} \leq NP/P \leq c_{2}$ in the infinity, corresponds to $NP \sim P$ (constant strength). For a hypoelliptic operator, we must assume $NP=P$ in the infinity implies $T=0$, that is $NP-P$ has locally isolated zeros, globally removable. Assume $P_{N}=\Pi N_{k}P$ and existence of $N_{k}P \equiv P$ on a segment, for $P$ not hypoelliptic. Thus, $P$ hypoelliptic implies for all $N \in (A)$, that $NP \neq P$, that is a continuation of $P$ does not contain any closed contours.
With the Fredholm operator representation, we have nontrivial kernel for $P$ not hypoelliptic and trivial
kernel for $P$ hypoelliptic (modulo regularizing action).

\newtheorem{prop5}[prop1]{Proposition}
\begin{prop5}
 Assume $\mu$ an analytic and locally reduced measure in $\mathcal{E}^{'(0)}$ such that the real and imaginary parts are orthogonal. Assume $\tilde{\mu}$, the
analytic continuation with algebraic segments and that the global representation has a trivial kernel. Then the corresponding functional is hypoelliptic.
\end{prop5}
(ref ArXiv)

\subsubsection*{Analytic continuations}
Assume that the dynamical system has right hand sides $(A,B)$ and $A=e^{v} P$ and $B=e^{v} Q$
so that $A/B=P/Q$. This gives an analytic continuation according to Lie. Assume $\frac{d}{d T} \gamma=(P,Q)$
and $\frac{d}{d T} \tilde{\gamma}=(A,B)$, such that $\frac{d}{d T} \tilde{x}= e^{v_{T}} \frac{d}{d T}x$,
with $\frac{d}{d T} v_{T}=0$ and $\frac{d}{d T} (e^{v_{T}} x)=e^{v_{T}} \frac{d}{d T}x$.

Assume now that $x^{\vartriangle}$ is a continuation of $x^{*}$, such that $x^{\vartriangle}=b(x^{*},y^{*})x^{*}$
and $y^{\vartriangle}=b(x^{*},y^{*})y^{*}$, where we assume $\{ b,b \} \neq 0$, under the relation
$x^{\vartriangle}.x + y^{\vartriangle}.y=b(x^{*}.x + y^{*}.y)$. We have that
$\frac{d y^{\vartriangle}}{d x^{\vartriangle}}=\frac{y^{*}}{x^{*}} (1 + \frac{\frac{d}{d t} \log y^{*}}{\frac{d}{d t} \log b})$.
Thus $\frac{y^{\vartriangle}}{x^{\vartriangle}}=\frac{y^{*}}{x^{*}}$ why if $x^{\vartriangle}+i y^{\vartriangle}$
is without sections, we have that $x^{*} + i y^{*}$ is without sections.

Consider for instance he system $(x,y) \rightarrow (X,Y) \rightarrow (M,W) \rightarrow (F_{1},F_{2})$ and consider continuation according to Lie $\frac{Y}{X} \rightarrow \frac{W}{M} \rightarrow \frac{F_{2}}{F_{1}} \sim \frac{1}{Q}$, where $Q$ is a polynomial and we assume the algebraicity is preserved. If $y'=g(y)=g(Qx)$ in the infinity, such that $g(Qx) \sim Q' g(x)$ in the infinity, where also $Q'$ is algebraic. The mapping $(x,y) \rightarrow x + i y$, assumes $x \bot y$, that is $\frac{y}{x} \sim \frac{1}{Q}$ in the infinity.
A sufficient condition is that $\frac{y}{x} \rightarrow 0$ in the infinity corresponding to algebraicity in the infinity. Orthogonality is written $<x,y>=0$ and onesidedness as $<x,y> \geq 0$.

\vsp

Consider the problem of finding a a factorization $\frac{\tilde{Y}}{\tilde{X}}=\frac{Q}{P}$,
where $P,Q$ are polynomials, such that $(P,Q) \rightarrow V$ is continuous. That is $P \tilde{Y}= Q \tilde{X}$. Note that given the continuation analytic, we have that 
$I(V \cup \Omega \backslash V)=I(V)I(\Omega \backslash V)$.  Note that presence of sections $Y= const X$
in the contraction condition, means symmetry $\frac{d G}{d x}=\frac{d G}{d y}$ which is implied by
$G(\ram x,y)=G(x, \ram y)$. Using the conditions $e^{-v} \frac{d \tilde{W}}{d x}=-Y + \frac{d v}{d x} e^{\phi}$
and $e^{-v} \frac{d \tilde{w}}{d y}=X + \frac{d v}{d y}e^{\phi}$, for $\tilde{W}=e^{\phi + v}$, we see
that $\frac{d v}{d x}=\frac{d v}{d y}=0$ implies a simple factorization that works. Further, symmetry for 
$\tilde{W}$ gives a contraction formula for the phase $Y/X \sim \frac{d v / d x}{d v / d y}$

\subsubsection*{Cousin's continuation method}

Cousin's method for analytic continuation can be used as follows. Let $\frac{d}{d T} (\tilde{F} - F) \sim \frac{d}{d T} \log \frac{\tilde{L}(\psi)}{<\eta_{T},\psi>} F=0$
iff $N(\psi)=<\eta_{T},\psi>=\tilde{L}(\psi)$ i $T=0$. 
More precisely, if we define an analytic $\alpha_{T}$, such that $\alpha_{T} \mid_{T=0} = <\eta,\psi>$ and $\alpha_{T} \mid_{T=\infty}=\tilde{L}(\psi)$ and such that the ramifier is regular and non-trivial. In this context, the FBI-transform is seen as analytic continuation over transversals of $f$.

\vsp

Assume the continuation has representation $f e^{\Sigma_{0}^{\infty} \phi_{j}}$, where we assume
$\frac{d}{d T} \phi_{j}=0$ and $\phi_{j} \neq 0$. Assume $V_{j}$ is the set where $\phi_{j} = 0$.
Then the sets $\cap_{j=0}^{N} V_{j}$ are continuation of support of $f$ with regular points.

Thus if $g_{j}$ gives group action $j=1,2$ such that $g_{1}=-g_{2}$ on $\Omega_{1} \cap \Omega_{2}$
then we have existence of $g$ with $g - g_{j} \in \mathcal{O}(\Omega_{j})$, thus given compatibility
conditions for sets of invariance, that is unique solution $g$ relative division of action.

\subsubsection*{Some remarks on Schur's lemma}

Note the lifting principle, where we have existence of $F$ holomorphic, such that $F(\gamma)(\zeta)=f(\zeta)$.
If $F$ is algebraic, then any spiral will be mapped onto a spiral.
In this case, we can form the mollification $F_{\delta} \rightarrow \delta_{0}$, and this implies $\zeta_{T} \rightarrow \zeta_{0}$
continuous. If we only assume $F \in H'$, we must consider $F(\gamma_{T})(\zeta) \rightarrow \delta_{0}$,
that is $<F(\gamma_{T}),\phi> \rightarrow < \delta,\phi>$. The algebraic geometry properties requires
analyticity, for instance $F=F_{1}F_{2}$ analytic and irreducible, then one factor, say $F_{1}$ is such that $F_{1} \sim \delta_{0}$.
If $F(\gamma)=F_{1}(\gamma)F_{2}(\gamma)$ implies $F_{1}(\gamma) \sim \delta_{0}$, then if $\gamma(0)=1$
$F_{1}(\gamma)(0)=1$. Assume now instead $H(\gamma)=F_{1}F_{2}(\gamma)-F_{1}(\gamma)F_{2}(\gamma)$.
Then $F_{1}^{-1}H(\gamma)=\big[ \delta,F_{2} \big](\gamma) - \delta_{0}(\gamma)F_{2}(\gamma)$ in $H'$
and if $\delta_{0}(\gamma)=1$, we have $\big[ \delta,F_{2} \big] \sim F_{2}$ why $F_{1}^{-1}H(\gamma)=0$.
In the same manner $H F_{2}^{-1}(\gamma)=\big[ F_{1},\delta \big](\gamma)-F_{1}(\gamma)\delta_{0}(\gamma)=0$
that is $HF_{2}=F_{1}H$. Using Schur's lemma, given $H$ linear we can determine $F_{1},F_{2}$.
The regularity properties for the corresponding inverses are determined by the involution condition.
Let $\tilde{\gamma}=F_{1}(\gamma)$ and $F(\gamma)=F_{1}(\tilde{\gamma})$. Assume $F_{2}(\gamma)$ corresponds to
polynomial right hand sides $(P,Q)$ and $F(\gamma)$ corresponds to analytic right hand sides $(X,Y)$
and that $P/Q \sim X/Y$ (Lie continuation). Then the lifting principle is solvable for $F_{1}$
but we can only guarantee existence $F$ in $H'$. If we assume an interpolation property, $P/Q \sim e^{\phi}P / e{^\phi}Y$
and $e^{\varphi} F(\gamma)=F_{1}(e^{\phi_{1}} \gamma)$, that is $e^{\varphi}Fe^{-\phi_{1}}$ is analytic, for $F \in H'$.
We can use Schur's lemma to determine $\phi,\varphi$.

\vsp

Assume $g \in G$ ($=\ram$) and let $< g f , \phi >=< f, {}^{t} g \phi>$. If $(I)=\{ f \quad g f = 0 \}$ implies $f({}^{t} g \phi)=0$ and $(J)=\{ \phi \quad \phi \bot g f \}$, where the ideals are considered in Schwartz type topology (\cite{Martineau}). If ${}^{t} g : (J) \rightarrow (I)$ is surjective, then $g$ is locally injective (and conversely). For instance if $i_{a}f=f+a$ is a continuous injection, then ${}^{t}i_{a}$ is surjective. Thus, $N({}^{t} g J)=\{ \zeta \quad {}^{t} g \phi(\zeta)=0 \quad \phi \in J \}=\{ \zeta_{T} \quad \phi(\zeta_{t}) \quad \phi \in J \}$. That is if we consider $g \in G$ as a homomorphism, we can form $(I)=(\mbox{ ker }g)$ and $\Omega=N(I)=\{ \zeta_{T} \quad \phi(\zeta_{T})=0 \}$

Note that using Schur's lemma, given two irreducible representations $\ram,r_{T}$ such that $\gamma r_{T}=\ram \gamma$, where $\gamma \neq 0$, we have that $\ram \sim r_{T}$. If we assume $\gamma$ absolute continuous, we have $\mid g \gamma - \gamma \mid < \epsilon$ iff $\mid \gamma({}^{t} g \zeta)-\gamma(\zeta) \mid < \epsilon'$, for positive $\epsilon,\epsilon'$ implies $\mid \zeta_{T} - \zeta \mid < \epsilon''$, if $g \in G$. Finally, we note that according to Schur's lemma if $\ram \gamma(\zeta)=\gamma(r_{T} \zeta)$
where $\zeta_{T},\gamma_{T}$ are complex lines. Then $\zeta_{T}$ can be determined uniquely
which means that the mapping $\gamma_{T}(\zeta) \rightarrow \zeta_{T}$ is continuous and locally 1-1.
Given the condition that all normal approximations are regular, we have that singularities are locally
isolated (the approximation property guarantees existence of normal approximations)

\subsubsection*{B\"acklund and the lifting principle}
Assume existence of $G$ such that $F(\widehat{M},\widehat{W})=G(\Top M,\Top W)$. The problem is to determine the regularity for $G$, given regularity for $F$.  Assume $d \sigma_{0}=(\mid X \mid^{2} + \mid Y \mid^{2})d x d y$ and $d \sigma_{1}=(\mid X \mid^{2}- \mid Y \mid^{2})d x d y$. Then $d \sigma_{0}=0$ implies $d \sigma_{1}=0$, then over $L^{1}(d \sigma_{1})$ according to Radon-Nikodym's theorem, we have existence of a Baire function $f_{0}$ such that 
$F(X,Y) d \sigma_{1}=F(X,Y)f_{0}(X,Y) d \sigma_{0}$. Let $p=(z,\overline{z})$, then the symmetry condition $G(z,\overline{z})=G(\overline{z},-z)$ is the condition $G(p)=G(\overline{p})$, which is a condition on the ramifier. Assume $\frac{d V}{d x}=-Y$ and $\frac{d V}{d y}=X$, then we have that on a set that is not reduced for contraction, $\mid \frac{d V}{d x} \mid=\mid \frac{d V}{d y} \mid$ implies $\mid Y \mid=\mid X \mid$. The condition $\int_{\Omega} F(X,Y) d \sigma_{0}=0$, for an algebraic $F$, means that $\sigma_{0}(\Omega)=0$ locally implies that $\Omega$ has zero measure. The condition $\int_{\Omega} F(X,Y) d \sigma_{1} \neq 0$, for an algebraic $F$ implies $\{ \mid X \mid=\mid Y \mid \} \subset \Omega^{c}$. Further, $\int_{\Omega} F d \sigma_{1}=0$, regardless of where $F$ is continuous. If $G(p)$ analytic with respect to $d \sigma_{0}$, then $G(p)$ is analytic outside $\{ \mid X \mid=\mid Y \mid \}$. Note that $(x+iy)(\overline{x}+i\overline{y})dx dy=d \sigma_{1}+i \mbox{ Re }(x \overline{y})dx dy$. Let $c(x,y)=x+iy$, then $c(\Top W)c(\overline{\Top W})=c(\Top M)c(\Top M^{\diamondsuit})$. We consider $(x,y) \rightarrow_{j_{1}} (\overline{x},\overline{y}) \rightarrow_{i_{2}} (\overline{x},-\overline{y})$ and $(x,y) \rightarrow_{i_{1}} (x,-y) \rightarrow_{j_{2}} (\overline{x},-\overline{y})$. Then $\overline{c(p)}=c(j_{1}i_{2}p)$ and $c(\overline{p})=c(j_{1}p)$. If the singular support for $G$ is one sided with respect to the $X-$ axes, then we may have that $G$ is analytic with respect to $\overline{c(p)}$ but not with respect to $c(\overline{p})$. Note that if $c(p)$ tracks a neighborhood of a circle, then $\log c(p)$ tracks a neighborhood of the real axes. If $X'=\Top M$ and $Y'=\Top W$ and $q=(X',Y')$, then $c(q)c(\overline{q}) \equiv 0$. Assume $M+iW \rightarrow_{\Top} \overline{\Top W}+i \Top W$ and $M +iW \rightarrow H+iG \rightarrow_{\rho} \overline{\Top W} + i \Top W$. We have for a Hamilton function $F$, with isolated singularities, using implicit derivation, existence of $\rho$ analytic, such that the domain of analyticity can be extended to $(\Top M,\Top W)$. The condition on isolated singularities for the dynamical system however requires that $F(\overline{X'},\overline{Y'})=F(\overline{X'},-\overline{Y'})$ does not hold. A sufficient condition is that $G > 0$. We note that the mapping $\Top$ can not be factorized into a part reduced with respect to orthogonality and a part reduced with respect to lineality. That is, there may be paths which appear in both sets simultaneously.

\subsection*{Continuation and group theory}

Lie's condition for linearity, $\frac{d^{2} \tilde{\phi}}{d x^{2}}=0$ implies $\tilde{\phi}$ is linear in $x$. Assume that $\frac{d^{2} \tilde{\phi}}{d x^{2}}=0$ iff $\frac{d \tilde{\phi}}{d x}=0$ in $L^{1}$. That is if for any $\tilde{\psi}$, we have existence of $\tilde{\phi}$, such that $\tilde{\psi}=\frac{d \tilde{\phi}}{d x}$, then we can write $M(\tilde{\psi}) \sim \tilde{\phi}$ in $L^{1}$, where $M$ denotes the arithmetic mean and we can assume $\tilde{\phi}$ has infinite order zero's (exponential representation). If $H(\tilde{\phi})=\frac{d}{d x} \delta * \tilde{\phi}=[I,\frac{d}{d x}](\tilde{\phi})=\frac{d}{d x}I(\tilde{\phi})$ (weak derivation). Assume $\Omega_{1}=\{ H(\tilde{\phi})=0 \}$ and $\Omega_{2}=\{ H^{2}(\tilde{\phi})=0 \}$. Then the condition $\frac{d^{2} \tilde{\phi}}{d x^{2}} + (\frac{d \tilde{\phi}}{d x})^{2}$ describes the set $\Omega_{1} \cap \Omega_{2} \neq \emptyset$. If we let $\tilde{\psi}=H(\tilde{\phi})$, we have $\frac{d}{d x} \tilde{\psi} + \tilde{\psi}^{2}=0$. Further, $MH(\tilde{\phi}) \sim I(\tilde{\phi})$ and the condition is $I(\tilde{\phi}) \sim -M_{2}H(\tilde{\phi})^{2}$ and
$\frac{H^{2}(\tilde{\phi})}{H(\tilde{\phi})}=-H(\tilde{\phi})$. Now let $V(\tilde{\phi})=\frac{H(\tilde{\phi})}{\tilde{\phi}}$, then $VH(\tilde{\phi})=-H(\tilde{\phi})$, that is $V$ corresponds on this set to reflection through the origin.

\vsp

If $T(\phi)(-\zeta)=-T(\phi)(\zeta)$ we say that $T$ is odd. Through change of
variable we have that $-T(\check{\phi}(\zeta)=T(\phi)(-\zeta)$ that is $T(\phi)=T(\check{\phi})$
implies $T$ odd. 
If we let $T(\check{\phi})=[T,I](\phi)$ and $[I,T](\phi)=-T(\phi)(-\zeta)$. Then we have
$[T,I](\phi)=[I,T](\phi)$ implies that $T$ is algebraic in $Exp$ (\cite{Martineau}). If now $H(\ram x,y)=H(x,\ram y)$ we have that $\int H(x,-y) \phi(y) d y= - \int H(x,y) \phi(-y) d y= - H(\check{\phi})(x)$,
thus $H(\phi)(-x)=-H(\check{\phi}(x)$. If on the other hand $H$ is not symmetric, we have that $-\int H(x,y) \phi(y) d y=\int H(x,y) \phi(y) d (-y)=\int H(x,-y) \phi(-y) d y= \check{H}(\check{\phi})(x)$

We note also the following example in representation theory, if $\Omega_{i}=N(I_{i})$ for $i=1,2$ and
$J=I_{1} + I_{2}$. Let $\mbox{ supp }I_{i}=\Omega_{i}^{c}=V_{i}$ for $i=1,2$ and assume $V_{i}$ are analytic sets.  
Then $\mbox{ supp }J=V_{1} \cup V_{2}$. We can write $V_{i}=N(e^{I_{i}^{0}})$ with $i=1,2$, and 
$N(e^{I_{1}^{0} + I_{2}^{0}})=V_{1} \cup V_{2}=\mbox{ supp }J$. Note that the representation is
dependent of the boundary condition (very regular).

Assume that $T=AB$ is the line with $T(\phi)=<ab,\phi>$. If $<T(\phi),\psi>=<<ab,\phi>,\psi>=<AB(\phi),\psi>=<\phi,{}^{t}B {}^{t}A \psi>$.
Thus if $T$ is linear, then ${}^{t}T$ linear. However ${}^{t} T \notin H'$. For instance if 
$FE(\phi)=\widehat{F}(\phi)$ and $F(e^{\phi})=e^{\tilde{F}\phi}$ and ${}^{t}E {}^{t}F(e^{\phi})={}^{t} \tilde{F}(\phi)=\log {}^{t} F E$.
Thus $\phi \in H$ implies $FE \in H'$ and ${}^{t}E {}^{t}F \notin H'$

\subsubsection*{Some remarks}

Two operators are characteristic in discussions about groups. If $c(x)=-x \in G$, for $x \in G$,
that is reflection through the origin is a group invariant property. The other property is given
$\phi_{\alpha} \in G$, this implies $\phi_{\alpha}^{\bot}=1-\phi_{\alpha} \in G$.

Assume $a = e^{\phi}$, $b = e^{\varphi}$ and $G \ni a,b$ implies $a b^{-1} \in G$. Then if $\log a,\log b \in G_{1}$, this
implies $\log a - \log b \in G_{1}$ and correspondingly for $G_{2} \ni \log \log a$.
The condition $N \log \phi \in G_{2}$ iff $\log \phi \in G_{2}$ implies $\phi^{N} \in G_{1}$ for all
$N$. Thus $\Sigma \phi^{N} \in G_{1}$ implies $(1-\phi)^{-1} \in G_{1} $ A group structure is said to be
isolated if $\phi^{N} \in G$ implies $\phi \in G$. If $a \bot b$ we have $(a+b)^{2}=a^{2} + b^{2}$.
The condition for Lagrange space is $\phi \in G$ iff $\phi^{\bot} \in G$. For Riesz groups we have that
the interpolation property implies $a < b_{1} + b_{2}$ which means existence of $a_{j} \leq b_{j}$,$j=1,2$
such that $a=a_{1} + a_{2}$. Thus, if $a \in G$ and $\phi=\log a \in L^{1}$ and $\phi \leq \mid \phi \mid$, we have $a^{\bot}=1-a \in G$.

The problem also is dependent on the following, if $\Sigma=\{ \varphi \equiv 0 \}$ analytic, we do not have
necessarily $\{ \varphi \equiv 1 \}$ analytic. Assume for this reason $\{ \varphi=d \varphi = 0 \}$ not isolated
points. If $\Sigma^{c}$ is analytic and $\{ \varphi \equiv 1 \}$ is closed in $\Sigma^{c}$, we can assume
that it is analytic. Particularly if $\Sigma$ is algebraic, we have that $\Sigma^{c}$ is locally analytic,
why $\overline{\{ \varphi \equiv 1 \}}$ is analytic. If $\varphi=\varphi_{1}$ on $\Sigma$ (algebraic)
means $\Sigma^{c}=\{ \varphi \neq \varphi_{1} \}$ is locally analytic, and $\overline{\{ \varphi - \varphi_{1} \equiv 1 \}}$ is
analytic. Thus, if $\Sigma,\Sigma^{c}$ is analytic, we have that we have existence of $\varphi_{1}$
analytic, such that $\varphi=1+\varphi_{1}$. Note that if $\Omega'=\{ \varphi \equiv 1 \}$ with $\frac{d}{d T} e^{\varphi} \equiv 0$
on $\Omega'$ and $\Omega=\{ \varphi \equiv 0 \}$ with $\frac{d}{d T} e^{\varphi} \equiv 0$ on $\Omega$,
then $\Omega'$ is not equivalent with $\Omega$. Further, if $G$ are transforms with underlying set $\Omega$, where
we assume $\Omega \ni x \rightarrow \frac{1}{x} \in \Omega$, implies $\Omega$ unbounded (functional representation),
that is $id(x)=x$, $id^{-1}(x)=\frac{1}{x}$ and $id^{-1} \in G$. 

We note the following problem, assume that $A,B,C$ are functionals with $\mbox{ supp }A \cap \mbox{ supp }B = \emptyset$
and that we have existence of $L=\mbox{ supp }C$ a line that separates $\mbox{ supp }A$ from $\mbox{ supp }B$.
Then we have $A \bot B$, $A \bot C$ and $B \bot C$.  If $A \bot B$, where 
$A$ has onesided support, in some sense also $B$ must have onesided support.
If $A$ is algebraic, there is a $E$ with $EA=I$ where $E$ can be assumed 
to have onesided support (with respect to a line $L$) If thus $B \bot A$ with $AB\varphi=0$ means $EAB \varphi=0$
which implies that $B$ has onesided support. Note that
if $A \bot B$ implies $A \otimes B=1$, then $(A \otimes B)^{2}=A^{2} \otimes B^{2}$. Further, if $A=A_{1}\ldots A_{N}$ and $B=B_{1}\ldots B_{N}$ and
we have existence of $A_{j}$ such that $A_{j} \bot B_{j}$ and such that $A \bot B$ as above. Thus if $A=A_{1}\ldots A_{N}$, where a factor is algebraic, then $A \bot B$ can be represented with onesided support.

\section{Sets of invariance}

Starting from the set of invariance associated to a group action, we discuss the continuation through
interpolation theory. Particularly if the distance is in time, there is not
an infinitesimal arbitrary decomposition of parameter intervals that corresponds to action. 
(``cold case'' ) why the integral representation must be motivated.

Assume $I$ is defined by $e^{<>}f=f \rightarrow \tau \widehat{f}=\widehat{f}$ and $J$ is defined
by $\tau f = f \rightarrow e^{<>} \widehat{f} = \widehat{f}$. Then, if we have that $N(I) \cap N(J)$
is removable $=N(I+J)$, this corresponds algebraically to a proposition on that the corresponding
homomorphy $g$ has symmetric kernel and $[I,g]=[g,I]$. 
Since $N(IJ)=N(I) \cup N(J)$, if $I \bot J$ then $(I+J)(I-J)=I^{2} + J^{2}$.
We should for this reason have that a decomposition in invariant sets is dependent of choice of 
boundary condition.

\subsubsection*{Lineality and the characteristic set}

Over the lineality, the class $(A)$ of continuations can be given by scalar products, that is $\{ e^{\Phi} \equiv 1 \}=\{ \Phi \equiv 0 \}$. If $\Phi=T(\psi)=0$ for all $\psi \in (J)$, then we have existence of $\eta$ such that $<\eta,\psi>=T(\psi)$, for instance $\eta=0$ or $\eta \in (J)^{\bot}$. A neighborhood of the lineality is given by $\{ \eta \quad <\eta,\psi> \geq 0 \}$. A bigger continuation can be given by $\{ \Phi \geq <\eta,\psi> \}. $
If $\log f + <\eta,\psi>=\log f$, then $<\eta,\psi>=0$ and conversely. Thus microlocal contribution can be seen as presence of a line in the infinity in $Z(<\eta,\psi>)$

\vsp

For $(I_{HE})$ with pseudo-base $F_{j}$ and $V=\cap N(F_{j})$, then there is a $\lambda$ such that if $f=0$ on $V$, we have $(\mbox{ Im }f)^{\lambda} \in (I_{HE})$. If we have presence of lineality and $\mbox{ Re }f \bot \mbox{ Im}f$, we must have $\mbox{ Im }f$ is bounded (that is we have unbounded sublevel sets). If $f$ is hypoelliptic, under the same conditions, we have that $\mbox{ Im }f$ is unbounded. 

\vsp

Assume $\Sigma=\{ f=c \}$, then a necessary condition for $\Sigma^{\bot}$ to be locally algebraic, is that  the leafs in $\Sigma$ are oriented, in the sense that $\gamma \bot \Sigma$ is locally on one side of the leaf in $\Sigma$. In this case we could write the lineality, $\Delta=\Delta_{+} \cup \Delta_{-}$, depending on what side of the leaf we are considering.

\vsp

Assume $p=(f,g) \rightarrow f+ig$, through $c(f,g)$, then $\overline{cp} \neq c \overline{p}$. We have $\int c(p)^{2} dz=0$, then $f^{2}=g^{2}$ and $\int c(p)c(\overline{p}) dz=0$ then $\mid f \mid=\mid g \mid$, if the integral is over a positive measured set. If $\int \mid c(p) \mid dz=0$, then $\mid f \mid=\mid g \mid=0$

Existence of lineality can be seen as a proposition of possibility to continue the symbol on a set of infinite order, that is the symbol is not reduced with respect to analytic continuation.  Assume $\rho(\Top \varphi)=\rho(\varphi^{*})$ on an algebraic set and $\rho(\Top \varphi)=0$ implies $\rho(\varphi^{*})=0$, then we have existence of $\varphi_{0}$ Baire such that $\rho(\Top \varphi)=\rho(\varphi_{0} \varphi^{*})$.

\subsubsection*{Some remarks}

Assume the representation $F_{T}=G_{T}+H_{T}$, where $G$ is reduced. For instance if we assume
$1/G_{T} \sim G_{1/T}$, we can assume $G_{T}$ reduced implies that it is locally algebraic in $1/T$. Assume ${}^{t}R_{T}=R_{T}$ so that $(HR_{T})^{t}=R_{T}{}^{t}H$ or better $H(\ram x,\ram y)=H(x,y)$, then $H$ can not be algebraic in $T$. Thus, if $F_{T}$ has lineality and the representation above with $G_{T}$ locally algebraic, then $H$ must be symmetric on an algebraic set.

\subsubsection*{Orthogonality}

Consider following interpolation property in $L^{1}$, $\frac{M(f^{(\lambda)})}{f^{(\lambda)}}=\rho_{\lambda}$, where $M$ is the arithmetic mean and
where $\frac{1}{f^{(\lambda)}} \sim \rho_{\lambda}\ldots \rho_{1}$. Thus $I \bot f^{(\lambda)}$ if 
$\exists \rho_{j} \rightarrow 0$, where the rest of $\rho_{j}$ are bounded. Note that if 
$\phi^{(\lambda)}=\log M(f^{(\lambda)})$ we have that $\frac{d}{dx} \phi^{(\lambda)} \sim \frac{1}{\rho_{\lambda}}$.
Note that if $f$ algebraic and $I \bot f^{(2)}$ we have that $I \bot f^{(1)}$. Further if $I \bot f$ 
and $\mid \mbox{ Im }f/ \mbox{ Re }f \mid < c$
then $I \bot \mbox{ Re }f$. If further on a disk $h \bot \mbox{ Im }f$ we have $h \bot \mbox{ Re }f$ 
and if $\mbox{ Im }f \bot \mbox{ Re }f$ in $L^{2}$, then $h \in C^{1}$ (\cite{AhlforsSario60}). 
The complement to $\sigma(f)$ (\cite{Mandelbrojt57}) is the set for representability as analytic function.
Precense of clustersets from lineality and orthogonality gives contribution to $\sigma(f)$
Assume the Riesz' representation theorem on $\sigma(f)$,
such that $<f,\varphi>=\mu_{n}$ over $\varphi \in C$, (moment) solvability means unique existence of $dV$ of bounded variation, such that 
$<f,\varphi>=\int \varphi d V$. Given an infinite sequence $\{ \mu_{n} \}$, the singularities that are given by $dV$, can be described by the convergence criteria associated to the problem.

\subsubsection*{The involution}
In connection with the conditions on involution, we note that $f(\zeta)=F(M,W)(\zeta)$ and $=F(M+iW)(\zeta)$. 
The condition $\mbox{ Im }F^{-1}f \bot M$ implies $< \mbox{ Im }f, \mbox{ Re }f>=0$ and 
$ \mbox{ Re }F^{-1}f \bot W$ implies $< \mbox{ Re }f,\mbox{ Im }f>=0$. For the first equality, we get on the 
right side $< \mbox{ Im }f, F(M,1)>=0$ or $F^{-1}(\mbox{ Re }f)=(M,1)$ and similarly for the other implication.

\vsp

Assume in $L^{1}$, $M(\varphi)=\varphi +v$, where here $M$ is the arithmetic mean, with $M(\varphi) \bot \psi$ and $\varphi \bot \psi$. Then we have
that $\frac{d}{d x}M(\varphi)=\frac{d}{d x} (\varphi + v) \sim \varphi$. Thus, $\frac{d v}{d x} \sim \varphi$ in $\infty$.
 If we assume $1 \bot \varphi$ and $\varphi \bot \frac{d}{d x} \varphi$,
if $\varphi'' / \varphi'$ bounded in $\infty$, then we have $\frac{d^{2} \varphi}{d x^{2}}/\varphi \rightarrow 0$
in $\infty$. Thus if $\varphi$ is algebraic in $\infty$ and $\varphi \bot \frac{d}{d x} \varphi$,
$1 \bot \varphi$, we have that $\varphi$ is reduced (HE). Consider the situation when $\varphi \bot \phi$,$\frac{d}{d x} \varphi \bot \phi$,
but $e^{-v_{2}} \frac{d^{2}}{d x^{2}} \varphi \bot \phi$, where $v_{2} > 0$. In this way we can construct a modified non algebraic reduced symbol.
\vsp

A transversal is a line that intersects two other lines. For instance $\varphi=c, \frac{d}{d x} \varphi=c$
(compare isolated singularities) that is $\frac{d}{d x}\varphi \bot \psi$ does not imply $\varphi \bot \psi$.
Thus for instance on both $\varphi'' \bot \phi$,$\varphi' \bot \phi$ does not imply $\varphi \bot \phi$
and $\varphi \bot \phi$,$\varphi' \bot \phi$ does not imply $\varphi'' \bot \phi$. The interpolation
property can now be compared with existence of $v$ such that $\varphi \bot e^{v} \phi$, where 
$v \geq 0$. If an infinitesimal movement is not planar, we can have $\varphi \bot \frac{d}{d x} \varphi,\varphi \bot \frac{d^{2}}{d x^{2}}\varphi$,
and so on up to infinite order. If $1 \bot \varphi$, $\varphi \bot \varphi'$ and $\varphi^{(j)}/\varphi^{(j-1)}$
bounded, for all $j$, for instance when $\varphi$ is a polynomial, we have that $\varphi^{(j)}/\varphi \rightarrow 0$,
for every $j$, that is $\varphi$ represents a reduced measure.

\vsp

Note that in the case with infinitely generated singularities, the class of surfaces that is generated
by normals is not equivalent with the one with surfaces generated by tangents. Further that a
oriented boundary is necessary for the normal to be locally algebraic.  

\subsubsection*{Examples}
Note that if
$\frac{d}{d T}f=(\frac{d}{d T}\phi) f$, so if $f \in L^{1}$ and 
$\frac{d}{d T} \phi$ is bounded, then $\frac{d}{d T}f \in L^{1}$. If 
$\mid \frac{d}{d T} \phi \mid \leq \mid \frac{d}{d T}f \mid$ and $f$ reduced, then if 
$\frac{d}{d T} f \in L^{1}$, we have $\frac{d}{d T} \phi \in L^{1}$. 

A condition necessary for continuation $N$ is $\log Nf$ is locally injective. This means that $\frac{\log N_{T}f}{ \log f}$ is constant implies $T=0$. Assume $f=e^{\phi}$ and $\frac{\log N_{T}f}{\log f}=\frac{<\eta,x> + \phi}{\phi}$. Thus, $\frac{<\eta_{T},x>}{\phi(x)}$ is constant implies $T=0$. The phase is defined where for instance $\phi(x) > <\eta_{T},x>$ locally. At the same time we assume there is a $\lambda_{T}$ such that $\phi(x) < <\lambda_{T},x>$.

Assume $\pi$ maps $\mid f \mid=1$ onto a planar domain, $(u,v) \in E$, where $f=u+iv$ and $u \bot v$.
Assume $\tilde{f}$ localized to the geometric set $V$ according to the reflection model, in such a manner that $\tilde{f}$ locally defined, is analytic over $V$, then the inverse mapping is locally defined, continuous. If $f(\zeta)=F(\gamma)(\zeta)$, note that if the problem is to localize the symbol to $V$, we only require that the continuations $\tilde{\gamma}(\zeta) \rightarrow \zeta \in V$
are continuous. But if we have to analyze $V$ from an analytic ideal, we require that $F(\tilde{\gamma})(\zeta)$ is locally analytic over $V$. If there is a difference of time involved,
then only the localization is possible, unless we can transport $F$ in time.

\newtheorem{main}[prop1]{Proposition}
\begin{main}
 Assume $\pi S^{1} \rightarrow E$, where $E$ is a planar domain and that $\varphi$ maps $E$ continuously onto $V$, a geometric set
of any order. Assume $f \in (I)(\Omega)$ a local geometric ideal and $\log \mid f \mid \in L^{1}(\Omega)$. If $\tilde{f}$ is $f$ localized to $V$ using the reflection model and $\pi$, where $\tilde{f}$ has TP between the reflection points, then the inverse mapping $\varphi$ is continuous $\tilde{f} \rightarrow V$.
\end{main}
(\cite{Lie96}),(\cite{Julia55}),(\cite{Collingwood66}),(\cite{BoutetdeMonvel76})

\bibliographystyle{amsplain}
\bibliography{ref}
\end{document}